\documentclass[10pt]{article}

\usepackage{amssymb}
\usepackage{eucal}
\usepackage{amsmath}
\usepackage{amsthm}
\usepackage{amscd}
\usepackage{graphics}
\usepackage{graphicx}
\usepackage{amsfonts}
\usepackage{latexsym}
\usepackage[all]{xy}

\numberwithin{equation}{section}
\newtheorem{thm}{\bf Theorem}[section]

\newtheorem{prop}[thm]{\bf Proposition}

\newtheorem{corollary}[thm]{Corollary}
\theoremstyle{remark}
\newtheorem{rem}{\bf Remark}[section]

\begin{document}

\title{Stability of equilibria for the $\mathfrak{so}(4)$ free rigid body}
\author{Petre Birtea$^a$, Ioan Ca\c{s}u$^b$, Tudor S. Ratiu$^c$, and Murat Turhan$^d$\\
{\small $^a$ Departamentul de Matematic\u a, Universitatea de Vest din Timi\c soara,}\\
{\small Bd. V Parvan, No. 4, RO--300223 Timi\c soara, Romania. E-mail: birtea@math.uvt.ro}\\
{\small $^b$ Departamentul de Matematic\u a, Universitatea de Vest din Timi\c soara,}\\
{\small Bd. V Parvan, No. 4, RO--300223 Timi\c soara, Romania.} \\
{\small Tel/fax: 0040-740-928382/0040-256-592316. E-mail: casu@math.uvt.ro (Corresponding Author)}\\
{\small $^c$ Section de Math\'ematiques and Bernoulli Center,}\\
{\small {\'E}cole Polytechnique F{\'e}d{\'e}rale de Lausanne, CH--1015 Lausanne, Switzerland. E-mail: tudor.ratiu@epfl.ch}\\
{\small $^d$ Section de Math\'ematiques,
{\'E}cole Polytechnique F{\'e}d{\'e}rale de Lausanne,}\\
{\small CH--1015 Lausanne,
Switzerland and}\\
{\small Yildiz Technical University,
Department of Mathematics,}\\
{\small Davutpasa Campus,
Esenler-Istanbul, Turkey.
E-mail: murat.turhan@epfl.ch}
}
\date{}

\maketitle



\begin{abstract}
The stability for all generic equilibria of the Lie-Poisson dynamics of the $\mathfrak{so}(4)$ rigid body
dynamics is completely determined.
It is shown that for the generalized rigid body certain Cartan
subalgebras (called of coordinate type) of
$\mathfrak{so}(n)$ are equilibrium points for the rigid body dynamics. In the case of
$\mathfrak{so}(4)$ there are three coordinate type Cartan subalgebras whose intersection with a regular adjoint orbit give three Weyl group orbits of equilibria. 
These coordinate type
Cartan subalgebras are the analogues of the three axes of equilibria
for the classical rigid body in $\mathfrak{so}(3)$.
In addition to these coordinate type Cartan equilibria there are others that come in curves.
\end{abstract}

\section{Introduction}

The goal of the present work is to find the analogue of the long axis--short axis stability theorem for the $ \mathfrak{so}(4)$-free rigid body. The first task is
to determine the analogue of the usual three axes of equilibria occuring in the dynamics of the classical 
$\mathfrak{so}(3) $ free rigid body.
It is shown that they are replaced by special Cartan subalgebras that we shall call \textit{coordinate type Cartan subalgebras}. For the general $\mathfrak{so}(n)$ rigid body it is proved that these coordinate type Cartan subalgebras are equilibria.

If $n = 4$, on a regular adjoint orbit, all the Cartan type equilibria are organized in three Weyl group orbits. 
The nonlinear stability and instability for these equilibria is determined taking into account the symplectic geometry of the orbit and the complete integrability of the system. The results in this paper complete and extend some previous work of Feh\'{e}r and Marshall \cite{FeMa03} and Spiegler
\cite{Spiegler04}. 

The standard energy methods for proving nonlinear stability (such as
energy-Casimir or energy-momentum) lead to very complicated computations and, in addition, for some equilibria the second variation of the augmented Hamiltonian is indefinite. Because of this, we take advantage of the symplectic geometry of the problem and its low dimensionality and we use Williamson's theorem (see e.g., \cite{BoFo04}), where stability is obtained by bringing the constants of motion in a normal form. We obtain three coordinate type Cartan subalgebras, $\mathfrak{t}_1$ , $\mathfrak{t}_2$, 
$\mathfrak{t}_3$ that, intersected with a regular adjoint orbit,  give rise to four equilibria $M^i_{\pm a,\pm b},M^i_{\pm b,\pm a}$, $i=1,2,3$. For a chosen ordering of the moments of inertia  we obtain the following stability results. The 
$\mathfrak{t}_2$-equilibria are unstable for the Lie-Poisson dynamics on $\mathfrak{so}(4)$ and unstable of center-saddle type for the dynamics on the generic adjoint orbit containing the equilibrium. The $\mathfrak{t}_3$-equilibria are stable for the Lie-Poisson dynamics on $\mathfrak{so}(4)$ and stable of center-center type on the generic adjoint orbit containing the equilibrium. In these two cases the Weyl group orbits of equilibria have the same stability or instability behavior. One would hope that this beautiful pattern holds in general. But this is not true, as bifurcations take place for the $\mathfrak{t}_1$-equilibria. The equilibria $M^1_{\pm a,\pm b}$ are stable for the Lie-Poisson dynamics on $\mathfrak{so}(4)$ and are stable of center-center type on the generic adjoint orbit containing 
$M^1_{\pm a,\pm b}$. The equilibria $M^1_{\pm b,\pm a}$ 
exhibit both stability and instability behavior that is
studied in detail in Section 5.

In addition to the Cartan type equilibria, there are, on every regular orbit, curves of equilibria. It is shown
that these are nonlinearly stable as a family, that is,
if a solution of the $\mathfrak{so}(4)$-free rigid body equation starts near an equilibrium on such a curve,
at any later time it will stay close to this curve but in the direction of the curve itself it may drift.

The implication of the topological structure of the energy-mementum level sets on bifurcation phenomena in the dynamics was extensively studied by Oshemkov \cite{Oshemkov87}, \cite{Oshemkov91}, Bolsinov and Fomenko \cite{BoFo04}.

\section{Equilibria for the generalized rigid body}

The equations of the rigid body on $\mathfrak{so}(n)$ are given by
\begin{equation}\label{rigideq}
\dot M=[M,\Omega], \end{equation} where $\Omega\in \mathfrak{so}(n)$,
$M=\Omega J+J\Omega\in \mathfrak{so}(n)$ with $J=\operatorname{diag}(\lambda_i)$, a real constant diagonal matrix satisfying $ \lambda _i + \lambda _j \geq 0 $, for all $i, j=1, \ldots, n$, $i \neq j $ (see, for example, \cite{Ratiu80}).
Note that $M=[m_{ij}]$ and $\Omega=[ \omega_{ij}]$ determine each other if and only if $\lambda_i + \lambda _j >0$ since $m_{ij} = ( \lambda_i+ \lambda_j) \omega_{ij}$ which physically means that the rigid body is not concentrated on a lower dimensional subspace of $\mathbb{R}^n$. Consequently, in the entire paper we shall assume that
 $\lambda_i+\lambda_j>0$ for all $i, j=1, \ldots, n$, $i \neq j $. In addition we shall study generic rigid bodies, i.e., all $\lambda_i$ are distinct.

It is well known and easy to verify that equations \eqref{rigideq} are Hamiltonian relative to the minus Lie-Poisson bracket 
\begin{equation}
\label{LP_n}
\{F,G\}(M) : = \frac{1}{2}\operatorname{Trace}(M [\nabla F(M), \nabla G(M)]),
\end{equation}
and the Hamiltonian function
\begin{equation}
\label{ham_n}
H(M) := - \frac{1}{4} \operatorname{Trace}(M \Omega).
\end{equation}
Here $F,G, H \in C ^{\infty}(\mathfrak{so}(n))$ and the gradient is taken relative to the Ad-invariant inner product
\begin{equation}
\label{ip_n}
\left\langle X, Y \right\rangle : = - \frac{1}{2}\operatorname{Trace}(XY), 
\qquad X, Y \in \mathfrak{so}(n)
\end{equation}
which identifies $(\mathfrak{so}(n))^*$ with $\mathfrak{so}(n)$. This means that $\dot{F} = \{F, H\}$ for all $F \in C^{\infty}(\mathfrak{so}(4))$, where $\{\cdot , \cdot \}$ is given by \eqref{LP_n} and $H $ by \eqref{ham_n}, if and only if \eqref{rigideq} holds.
Note that the linear isomorphism $X \in \mathfrak{so}(n) \mapsto XJ + JX 
\in \mathfrak{so}(n)$ is self-adjoint relative to the inner product \eqref{ip_n} and thus $\nabla H(M) =  \Omega$.

Let $E_{ij}$ be the constant antisymmetric matrix with $1$ on line
$i$ and column $j$ when $i<j$, that is, the $(k,l) $-entry of $ E_{ij}$ equals $(E_{ij})_{kl} = \delta_{ki} \delta_{lj} - \delta_{kj} \delta_{li}$. Then $\{E_{ij} \mid i<j\}$ is a basis for the Lie algebra $\mathfrak{so}(n)$. We have
\[
\left(E_{ij}E_{ks} \right)_{ab} = \delta_{ai} \delta_{jk} \delta_{bs}
- \delta_{aj} \delta_{ik} \delta_{bs} 
- \delta_{ai} \delta_{js} \delta_{bk}
+ \delta_{aj} \delta_{is} \delta_{bk}
\] 
and hence $E_{ij}^2$ is the diagonal matrix whose only non-zero entries $-1$ occur on the $i$th and $j$th place. In addition, if $i<j$ and $k<s$, we get
\begin{equation}
\label{bracket}
\left[E_{ij}, E_{ks} \right] = \delta_{jk}E_{is} + \delta_{is}E_{jk} - \delta_{ik} E_{js} - \delta_{js} E_{ik}
\end{equation}
where $E_{rp} : = -E_{pr}$, if $r>p$.

Since 
\[
[M,\Omega]=[\Omega J+J\Omega,\Omega]
=\Omega J\Omega+J\Omega^2-\Omega^2 J-\Omega J\Omega
=[J,\Omega^2].
\]
we see that $M $ is an equilibrium if and only if $[J,\Omega^2]=0$. 

Since we assume that all $\lambda_i $ are distinct, this condition is equivalent to the statement that  $\Omega^2$ is a diagonal matrix.

\begin{thm} \label{C equilibria}
Let $\mathfrak{h}\subset \mathfrak{so}(n)$ be a
Cartan subalgebra whose basis is a subset of $\{E_{ij} \mid i<j\}$. Then any
element of $\mathfrak{h}$ is an equilibrium point of the rigid body equations \eqref{rigideq}.
\end{thm}

\noindent{\bf Proof.} We have to prove that for any $M\in\mathfrak{h}$ the
matrix $\Omega^2$ is diagonal. Since the linear isomorphism $\Omega \leftrightarrow M$ is given by a diagonal matrix in the basis 
$\{E_{ij} \mid i<j\}$ of $\mathfrak{so}(n)$ it follows that $M \in \mathfrak{h}$ if and only if $\Omega \in \mathfrak{h}$.

So let $\Omega\in\mathfrak{h}$ with
$\Omega=\sum\limits_{s=1}^k\alpha_sE_s$, where $k:=[n/2] = \dim \mathfrak{h}$ and
$\{E_1,...,E_k\} \subset  \{E_{ij} \mid i<j\}$ is the basis of $\mathfrak{h}$.
 Then
\begin{equation}
\label{rb_comp}
\Omega^2=\left(
\sum\limits_{s=1}^k \alpha_sE_s\right)^2
=\sum\limits_{s=1}^k\alpha_s^2E_s^2+\sum\limits_{l\not=
p}\alpha_l\alpha_p(E_lE_p+E_pE_l).
\end{equation}
Since $\mathfrak{h}$ is a Cartan
subalgebra we have $[E_l,E_p]=0$ which is equivalent to
$E_lE_p=E_pE_l$ for any $l,p \in  \{1,\ldots, k\}$. Then
$(E_lE_p)^t=E_p^tE_l^t=(-1)^2E_pE_l=E_lE_p$. Consequently, the
matrix $E_lE_p$ is symmetric. Since $E_l, E_p \in  \{E_{ij} \mid i<j\}$, we distinguish the following cases for $l \neq p$:

(a) $E_l=E_{ij},E_p=E_{js},i<j<s$, in which case the product $ E_lE_p$ is not symmetric because the $ (i,s) $-entry equals $ 1 $ and the $ (s,i) $-entry vanishes. So this case cannot occur.

(b) $E_l=E_{ij},E_p=E_{sj},i<j,s<j,i\not= s$. Then, $ E_{ij}E_{sj}$ is not symmetric because the $ (i,s) $-entry equals $-1 $ and the $ (s,i) $-entry vanishes. So this case cannot occur.

(c) $E_l=E_{ij},E_p=E_{is},i<j,i<s,j\not= s$. Then $ E_{ij}E_{is}$ is not symmetric because the $ (j,s)$-entry equals $ -1 $ and the $ (s,j)$-entry vanishes. So this case cannot occur.

(d) $E_l=E_{ij},E_p=E_{ks},i<j,k<s,\{i,j\}\cap \{k,s\}=\varnothing$.
In this case $E_lE_p=O_n$.

Thus, the only possible case in \eqref {rb_comp} is (d) which implies that
$$\Omega^2=\sum\limits_{s=1}^k\alpha_s^2E_s^2$$
which is a diagonal matrix. \rule{0.5em}{0.5em}
\medskip

We shall call a Cartan subalgebra as in Theorem
\ref{C equilibria} a \textit{coordinate type Cartan subalgebra}.
The dynamics of (\ref{rigideq}) leaves the adjoint
orbits of $\operatorname{SO}(n)$ invariant. Since the intersection of a regular orbit (that is, one passing through a regular semisimple element of $ \mathfrak{so}(n)$) with a Cartan subalgebra is a Weyl group orbit (see, e.g. \cite{Kostant73}), we conclude that the union of the Weyl group orbits determined by the coordinate type Cartan subalgebras of $ \mathfrak{so}(n)$ are equilibria of \eqref{rigideq}. As we shall see,  if $n\geq 4$, the system \eqref{rigideq} has also other equilibria that are not coming from coordinate type Cartan subalgebras.

\section{The adjoint orbits of $\mathfrak{so}(4)$}

The Lie algebra of the compact subgroup 
$\operatorname{SO}(4)=\{A\in \mathfrak{gl}(4, \mathbb{R}) 
\mid A^tA=I_4,\det (A)=1\}$ of the special linear Lie group
$\operatorname{SL}(4,\mathbb{R})$ is $\mathfrak{so}(4)$.
In this section we present the geometry of the (co)adjoint
orbits of $\operatorname{SO}(4)$ in $\mathfrak{so}(4)$.

We choose as basis of $\mathfrak{so}(4)$ the matrices
$E_1=-E_{23}$, $E_{2}=E_{13}$, $E_{3}=-E_{12}$, $E_{4}=E_{14}$, $E_{5}=E_{24}$, $E_{6}=E_{34}$
and hence we represent $\mathfrak{so}(4)$ as
\begin{equation}
\label{so4_representation}
\mathfrak{so}(4)=\left\{\left.
M=\left[\begin{array}{cccc}
0&-x_3&x_2&y_1\\
x_3&0&-x_1&y_2\\
-x_2&x_1&0&y_3\\
-y_1&-y_2&-y_3&0\end{array}\right] \, \right|\, x_1,x_2,x_3,y_1,y_2,y_3\in \mathbb{R}
\right\}.
\end{equation}
Since $\hbox{rank}\,\mathfrak{so}(4) = 2$, there are two functionally independent Casimir functions for the minus Lie-Poisson bracket, which are given by
\[
C_1(M) := -\frac{1}{4}\operatorname{Trace}(M^2)
= \frac{1}{2}\left(\sum\limits\limits_{i=1}^3 x_i^2+\sum\limits\limits_{i=1}^3 y_i^2 \right)
\]
and
\[
C_2(M) := - \operatorname{Pf}(M)
=\sum\limits\limits_{i=1}^3 x_iy_i.
\]
Thus the generic adjoint orbits are the level sets
$$
\operatorname{Orb}_{c_1c_2}(M)=(C_1\times C_2)^{-1}(c_1,c_2), 
\qquad (c_1,c_2) \in \mathbb{R}^2.
$$
Note that if $M \neq 0$, then $\mathbf{d}C_j(M)  \neq 0$  for $j=1,2$.

The Lie algebra $\mathfrak{so}(4) = \mathfrak{so}(3) \times \mathfrak{so}(3)$ is
of type $A_{1}\times A_{1}$ and, consequently, the positive Weyl
chamber, which is the moduli space of (co)adjoint orbits,
is isomorphic to the positive quadrant in $\mathbb{R}^2$. In the
basis of $\mathfrak{so}(4)$ that we have chosen above, the positive Weyl
chamber is given by the set
$\{(c_1,c_2)\in \mathbb{R}^2\mid c_1\geq |c_2| \}$.

To characterize the adjoint orbits of $\operatorname{SO}(4)$ it is convenient to split $\mathfrak{so}(4) = V_{ \mathbf{u}} \oplus V_{ \mathbf{v}}$, where $V_{ \mathbf{u}}: = \{\mathbf{u}: =(u_1, u_2, u_3)\mid u_i:= x_i+y_i, i=1,2,3\}\cong \mathbb{R}^3$, $V_{ \mathbf{v}}: = \{\mathbf{v}: = (v_1, v_2, v_3)\mid v_i:= x_i-y_i, i=1,2,3\} \cong \mathbb{R}^3$. Instead of the independent Casimir functions $C_1$, $C_2$ we consider the following two independent Casimir functions
\[
D_{\mathbf{u}}(M) : = 2C_1(M)+2C_2(M) = \| \mathbf{u}\|^2, \qquad 
D_{\mathbf{v}}(M) : = 2C_1(M)-2C_2(M) = \|\mathbf{v} \|^2. 
\]
Note that
\[
\operatorname{Orb}_{c_1c_2}(M)=(C_1\times C_2)^{-1}(c_1,c_2)
= (D_{ \mathbf{u}} \times D_{ \mathbf{v}}) ^{-1}(d_1, d_2), \quad \text{where} \quad d_1 = 2c_1+2c_2,\quad d_2 = 2c_1-2c_2.
\]

These considerations yield the following characterization of the  $\operatorname{SO}(4)$-adjoint orbits. 

\begin{thm} \label{caracterizare}
Denote by $S^2_r$ the sphere in $\mathbb{R}^3$ of radius $r $.
If $c_1>0$ and $c_1>|c_2|$, then the adjoint orbit
$\operatorname{Orb}_{c_1c_2}(M)$ equals
$S^2_{\sqrt{2c_1+2c_2}}\times S^2_{\sqrt{2c_1-2c_2}}$, where $S^2_{\sqrt{2c_1+2c_2}} \subset V_{ \mathbf{u}}$, $S^2_{\sqrt{2c_1-2c_2}} \subset V_{ \mathbf{v}}$, and hence it is regular. If $c_1=|c_2|>0$, then
the adjoint orbit $\operatorname{Orb}_{c_1c_2}(M)$ is
either $S^2_{2\sqrt{c_1}} \times \{0\}$, with $S^2_{2\sqrt{c_1}} \subset V_{ \mathbf{u}}$, or $\{0\} \times S^2_{2\sqrt{c_1}}$, with $S^2_{2\sqrt{c_1}} \subset V_{ \mathbf{v}}$, 
and so it is singular. If $c_1=c_2=0$,
then the adjoint orbit $\operatorname{Orb}_{c_1c_2}$ is
the origin of $\mathfrak{so}(4)$ and so it is singular.
\end{thm}

In all that follows we shall denote by $\operatorname{Orb}_{c_1;c_2}$ the regular adjoint orbit $\operatorname{Orb}_{c_1c_2}$,
where $c_1>0$ and $c_1>|c_2|$, which is equivalent to $d _1, d _2>0$.

A straightforward computation shows that the coordinate type Cartan subalgebras of $\mathfrak{so}(4)$ in the chosen basis are
$\mathfrak{t}_1,\mathfrak{t}_2,\mathfrak{t}_3$, where
$$\mathfrak{t}_1:=\operatorname{span}(E_1,E_4)
=\left\{\left. M_{a,b}^1:=\left[\begin{array}{cccc}
0&0&0&b\\
0&0&-a&0\\
0&a&0&0\\
-b&0&0&0\end{array}\right]\, \right| \,a,b\in \mathbb{R}
 \right\},$$
$$\mathfrak{t}_2:=\operatorname{span}(E_2,E_5)
=\left\{\left.  M_{a,b}^2:=\left[\begin{array}{cccc}
0&0&a&0\\
0&0&0&b\\
-a&0&0&0\\
0&-b&0&0\end{array}\right] \, \right| \, a,b\in \mathbb{R}
 \right\},$$
$$\mathfrak{t}_3:=\operatorname{span}(E_3,E_6)
=\left\{\left. M_{a,b}^3:=\left[\begin{array}{cccc}
0&-a&0&0\\
a&0&0&0\\
0&0&0&b\\
0&0&-b&0\end{array}\right] \, \right| \,a,b\in \mathbb{R}
 \right\}.$$

The intersection of a regular adjoint orbit and a coordinate
type Cartan subalgebra has four elements which represents a Weyl group orbit. Thus we expect at least twelve equilibria for the rigid body equations \eqref{rigideq} in the case of $\mathfrak{so}(4) $. Specifically, we have the following result.

\begin{thm} \label{Weil_group_orbit}
The following equalities hold:
\begin{itemize}
\item[{\rm (i)}] $\mathfrak{t}_1\cap \operatorname{Orb}_{c_1;c_2}=\left\{
M_{a,b}^1,M_{-a,-b}^1,M_{b,a}^1,M_{-b,-a}^1
\right\}$,
\item[{\rm (ii)}]  $\mathfrak{t}_2\cap \operatorname{Orb}_{c_1;c_2}=\left\{
M_{a,b}^2,M_{-a,-b}^2,M_{b,a}^2,M_{-b,-a}^2
\right\}$,
\item[{\rm (iii)}]  $\mathfrak{t}_3\cap
    \operatorname{Orb}_{c_1;c_2}=\left\{
    M_{a,b}^3,M_{-a,-b}^3,M_{b,a}^3,M_{-b,-a}^3
    \right\}$,
\end{itemize}
where
\begin{equation}
\label{values_a_b}
\left\{
\begin{aligned}
a&=\frac{1}{\sqrt{2}}\left(\sqrt{c_1+c_2}+\sqrt{c_1-c_2}\right) \\
b&=\frac{1}{\sqrt{2}}\left(\sqrt{c_1+c_2}-\sqrt{c_1-c_2}\right).
\end{aligned}
\right.
\end{equation}
\end{thm}

\noindent {\bf Proof.} 
Let $M_{\alpha,\beta}^1\in \mathfrak{t}_1\cap \operatorname{Orb}_{c_1;c_2}$.
Then $M_{\alpha,\beta}^1\in \mathfrak{t}_1$, $2c_1 = 2C_1(M_{\alpha,\beta}^1)=\alpha^2 + \beta^2$,
and $c_2 = C_2(M_{\alpha,\beta}^1)=\alpha\beta$. This system of equations has the solutions
$$(\alpha,\beta)\in\{(a,b),(-a,-b),(b,a),(-b-a)\},$$
where
$$\left\{
\begin{aligned}
a&=\frac{1}{\sqrt{2}}\left(\sqrt{c_1+c_2}+\sqrt{c_1-c_2}\right) \\
b&=\frac{1}{\sqrt{2}}\left(\sqrt{c_1+c_2}-\sqrt{c_1-c_2}\right).
\end{aligned}
\right.$$
Similar arguments with obvious modifications prove assertions
(ii) and (iii). \rule{0.5em}{0.5em}
\medskip

The intersections $\mathfrak{t}_1\cap \operatorname{Orb}_{c_1;c_2}$,  
$\mathfrak{t}_2\cap \operatorname{Orb}_{c_1;c_2}$, 
$\mathfrak{t}_3\cap \operatorname{Orb}_{c_1;c_2}$ are Weyl group orbits.

\section{The $\mathfrak{so}(4)$-rigid body}

We recall that we work under the generic assumptions $\lambda_i+ \lambda_j>0$ for $i \neq j$ and all $\lambda_i $ are distinct. The equations of motion are $\dot{M}=[M,\Omega]$, where $M=J\Omega+\Omega J$, $\Omega \in  \mathfrak{so}(4)$, $J=\operatorname{diag}\{\lambda_1,\lambda_2,\lambda_3,\lambda_4\}$,
$\lambda_1,\lambda_2,\lambda_3,\lambda_4\in \mathbb{R}$.
The relationship between $\Omega=[\omega_{ij}]\in \mathfrak{so}(4)$ and the matrix $M\in \mathfrak{so}(4)$ in the representation \eqref{so4_representation} is hence given by
$$\left.\begin{array}{lll}
(\lambda_3+\lambda_2)\omega_{32}=x_1& \quad 
(\lambda_1+\lambda_3)\omega_{13}=x_2& \quad
(\lambda_2+\lambda_1)\omega_{21}=x_3\\
(\lambda_1+\lambda_4)\omega_{14}=y_1& \quad 
(\lambda_2+\lambda_4)\omega_{24}=y_2& \quad 
(\lambda_3+\lambda_4)\omega_{34}=y_3
\end{array}\right.$$
and thus the equations of motion \eqref{rigideq} are equivalent for $n=4$ to the system
\begin{equation}\label{soM}
\left\{\begin{array}{l}
\vspace{.1cm}
\dot x_1=\left(\frac{1}{\lambda_1+\lambda_2}-\frac{1}{\lambda_1+\lambda_3}\right)x_2x_3+\left(\frac{1}{\lambda_3+\lambda_4}-\frac{1}{\lambda_2+\lambda_4}\right)y_2y_3\\
\vspace{.1cm}
\dot x_2=\left(\frac{1}{\lambda_2+\lambda_3}-\frac{1}{\lambda_1+\lambda_2}\right)x_1x_3+\left(\frac{1}{\lambda_1+\lambda_4}-\frac{1}{\lambda_3+\lambda_4}\right)y_1y_3\\
\vspace{.1cm}
\dot x_3=\left(\frac{1}{\lambda_1+\lambda_3}-\frac{1}{\lambda_2+\lambda_3}\right)x_1x_2+\left(\frac{1}{\lambda_2+\lambda_4}-\frac{1}{\lambda_1+\lambda_4}\right)y_1y_2\\
\vspace{.1cm}
\dot y_1=\left(\frac{1}{\lambda_3+\lambda_4}-\frac{1}{\lambda_1+\lambda_3}\right)x_2y_3+\left(\frac{1}{\lambda_1+\lambda_2}-\frac{1}{\lambda_2+\lambda_4}\right)x_3y_2\\
\vspace{.1cm}
\dot y_2=\left(\frac{1}{\lambda_2+\lambda_3}-\frac{1}{\lambda_3+\lambda_4}\right)x_1y_3+\left(\frac{1}{\lambda_1+\lambda_4}-\frac{1}{\lambda_1+\lambda_2}\right)x_3y_1\\
\vspace{.1cm}
\dot y_3=\left(\frac{1}{\lambda_2+\lambda_4}-\frac{1}{\lambda_2+\lambda_3}\right)x_1y_2+\left(\frac{1}{\lambda_1+\lambda_3}-\frac{1}{\lambda_1+\lambda_4}\right)x_2y_1.
\end{array}\right.
\end{equation}
The Hamiltonian \eqref{ham_n} has in the form
\begin{align*}
H(M)&=-\frac{1}{4}\hbox{Trace}(M\Omega)\\
&=\frac{1}{2}\left(\frac{1}{\lambda_2+\lambda_3}x_1^2+\frac{1}{\lambda_1+\lambda_3}x_2^2+\frac{1}{\lambda_1+\lambda_2}x_3^2+\frac{1}{\lambda_1+\lambda_4}y_1^2+\frac{1}{\lambda_2+\lambda_4}y_2^2+\frac{1}{\lambda_3+\lambda_4}y_3^2\right).
\end{align*}
The following result gives the list of all equilibria of the system \eqref{soM}. The proof is a straightforward computation.
\begin{thm} \label{echilibre} 
If $\mathcal{E}$ denotes the set of the equilibrium points
of \eqref{soM}, then
${\cal E}=\mathfrak{t}_1\cup \mathfrak{t}_2\cup \mathfrak{t}_3\cup \mathfrak{s}_+ \cup \mathfrak{s}_-$,
where $\mathfrak{s}_{\pm}$ are the three dimensional vector subspaces given by
\begin{equation*}
\mathfrak{s}_{\pm}:=\operatorname{span}_{\mathbb{R}}\left\{\left(\frac{1}{\lambda_1+\lambda_4}E_1 \pm \frac{1}{\lambda_2+\lambda_3}E_4\right),\,\left(\frac{1}{\lambda_2+\lambda_4}E_2\pm \frac{1}{\lambda_1+\lambda_3}E_5\right), \,\left(\frac{1}{\lambda_3+\lambda_4}E_3\pm \frac{1}{\lambda_1+\lambda_2}E_6\right) \right\}.
\end{equation*}
\end{thm}

\medskip

Note that $\mathfrak{s}_{\pm}$ are not Lie subalgebras of $\mathfrak{so}(4)$ and that $\mathfrak{s}_{\pm}\cap \mathfrak{t} _i \neq\varnothing$, $i=1,2,3$. Let us compute $\mathfrak{s}_{\pm}\cap \operatorname{Orb}_{c_1c_2}$. 
For an arbitrary element of $\mathfrak{s}_{\pm}$, we have
\begin{equation}
\label{coord_s}
\begin{array}{ccc}
\displaystyle
x_1 = \frac{a_1}{ \lambda_1+ \lambda_4}, &
\displaystyle
\quad x_2 = \frac{a_2}{ \lambda_2+ \lambda_4}, &
\displaystyle
\quad x_3 = \frac{a_3}{ \lambda_3+ \lambda_4},\\
\\
\displaystyle
 y_1 = \pm \frac{a_1}{ \lambda_2+ \lambda_3}, & 
 \displaystyle 
 \quad  y_2 = \pm \frac{a_2}{ \lambda_1+ \lambda_3}, &
 \displaystyle
 \quad y_3 = \pm \frac{a_3}{\lambda_1+ \lambda_2}.
\end{array} 
\end{equation}
Note that if $a _i = a _j = 0 $, $i \neq j$, $i, j \in \{1,2,3\}$, then the equilibrium lies in $\mathfrak{t}_k $, where $k \in  \{1,2,3\} \setminus \{i,j\}$. Thus the equilibria in $\mathfrak{s}_{\pm}$ that are not in $\mathfrak{t}_1\cup \mathfrak{t}_2\cup \mathfrak{t}_3$ must have at least two of $a_1, a _2, a _3$ different from zero.
From \eqref{coord_s} we deduce
\begin{align}
\label{first_condition}
0< c_1&=C_1(M) = 
\frac{1}{2}\left(\sum\limits\limits_{i=1}^3 x_i^2
+\sum\limits\limits_{i=1}^3 y_i^2 \right) \nonumber \\
&= \frac{a_1^2 }{2}\left(\frac{1}{(\lambda_1+ \lambda_4)^2} + \frac{1}{(\lambda_2 + \lambda_3)^2} \right) 
+ \frac{a_2^2}{2} \left(\frac{1}{(\lambda_2+ \lambda_4)^2} + \frac{1}{ (\lambda_1 + \lambda_3)^2} \right) \nonumber \\
& \qquad \qquad
+ \frac{a_3^2}{2} \left(\frac{1}{(\lambda_3+ \lambda_4)^2} + \frac{1}{( \lambda_1 + \lambda_2)^2} \right) 
 \end{align}
and
\begin{equation}
\label{second_condition}
c_2=C_2(M) 
=\sum\limits\limits_{i=1}^3 x_iy_i
=\pm \left(\frac{a_1^2}{(\lambda_1+ \lambda_4)( \lambda_2 + \lambda_3)} + \frac{a_2^2}{(\lambda_2+ \lambda_4)( \lambda_1+ \lambda_3)} + \frac{a_3^2}{(\lambda_3+ \lambda_4)( \lambda_1+ \lambda_2)} \right)
\end{equation}
which shows that $c_1 \geq |c_2|> 0$.  Since we are working in the
generic case when all $\lambda_i$ are distinct, if follows that at least one 
of the conditions 
\begin{equation}
\label{genericity_conditions}
\lambda_1+ \lambda_4 \neq \lambda_2 + \lambda_3,  \qquad \lambda_2+ \lambda_4 \neq \lambda_1 + \lambda_3, \qquad \lambda_3+ \lambda_4 \neq \lambda_1 + \lambda_2
\end{equation}
 holds and hence we obtain $c_1 > |c_2|>0$ which shows that  all equilibria in $\mathfrak{s} _{\pm}$ that are not in $\mathfrak{t}_1 \cup \mathfrak{t}_2 \cup \mathfrak{t}_3$ necessarily lie on a regular adjoint orbit 
$\operatorname{Orb}_{c_1; c_2}$ (see Theorem \ref{caracterizare}).   
These equilibria are not isolated. To describe them, express conditions \eqref{first_condition} and \eqref{second_condition} in the new variables
\[
b_1: = \frac{a_1}{\sqrt{( \lambda_1+ \lambda_4)( \lambda_2+ \lambda_3)}}\,, \qquad 
b_2: = \frac{a_2}{\sqrt{( \lambda_2+ \lambda_4)( \lambda_1+ \lambda_3)}}\,, \qquad 
b_3: = \frac{a_3}{\sqrt{( \lambda_3+ \lambda_4)( \lambda_1+ \lambda_2)}}\,,
\]
and hence $\operatorname{Orb}_{c_1; c_2} \cap \mathfrak{s}_{\pm}$ are
intersections of ellipsoids with spheres given by
\begin{align}
0<c_1 
&= \frac{b_1^2 }{2}\left(\frac{\lambda_2+ \lambda_3}{\lambda_1+ \lambda_4} + \frac{\lambda_1+ \lambda_4}{\lambda_2 + \lambda_3} \right) 
+ \frac{b_2^2}{2} \left(\frac{\lambda_1+ \lambda_3}{\lambda_2+ \lambda_4} + \frac{\lambda_2+ \lambda_4}{\lambda_1 + \lambda_3} \right) + \frac{b_3^2}{2} \left(\frac{\lambda_1+ \lambda_2}{\lambda_3+ \lambda_4} + \frac{\lambda_3+ \lambda_4}{\lambda_1 + \lambda_2} \right)
\label{first_condition_b}\\
c_2 &= \pm\left(b_1^2+ b_2 ^2+ b_3 ^3\right).
\label{second_condition_b}  
\end{align}
This also shows that if $c_2>0 $ then $\mathfrak{s} _- \cap \operatorname{Orb}_{c_1; c_2} = \varnothing $ and that if $c_2<0 $, then 
$\mathfrak{s}_+ \cap \operatorname{Orb}_{c_1; c_2} = \varnothing $. 
The considerations above prove the following result, which is depicted by Fig. 1.

\begin{corollary} \label{cor_equ_orbit}
On a generic adjoint orbit $\operatorname{Orb}_{c_1;c_2}$, $c_1>|c_2|>0$, the equilibria of \eqref{soM} are given by the twelve points in Theorem \ref{Weil_group_orbit} forming three Weyl group orbits in $\mathfrak{t}_1$, $\mathfrak{t}_2$, $\mathfrak{t}_3$, and the subsets in $\mathfrak{s}_{\pm}$ described by \eqref{first_condition_b} and \eqref{second_condition_b}. On a generic adjoint orbit $\operatorname{Orb}_{c_1;c_2}$, $c_1>|c_2|=0$, we have as equilibria for \eqref{soM} only the twelve isolated points in Theorem \ref{Weil_group_orbit} forming three Weil group orbits in $\mathfrak{t}_1$, $\mathfrak{t}_2$, 
$\mathfrak{t}_3$.
\end{corollary}

\begin{figure}[!ht]
\begin{center}
\includegraphics[scale=0.3,angle=0]{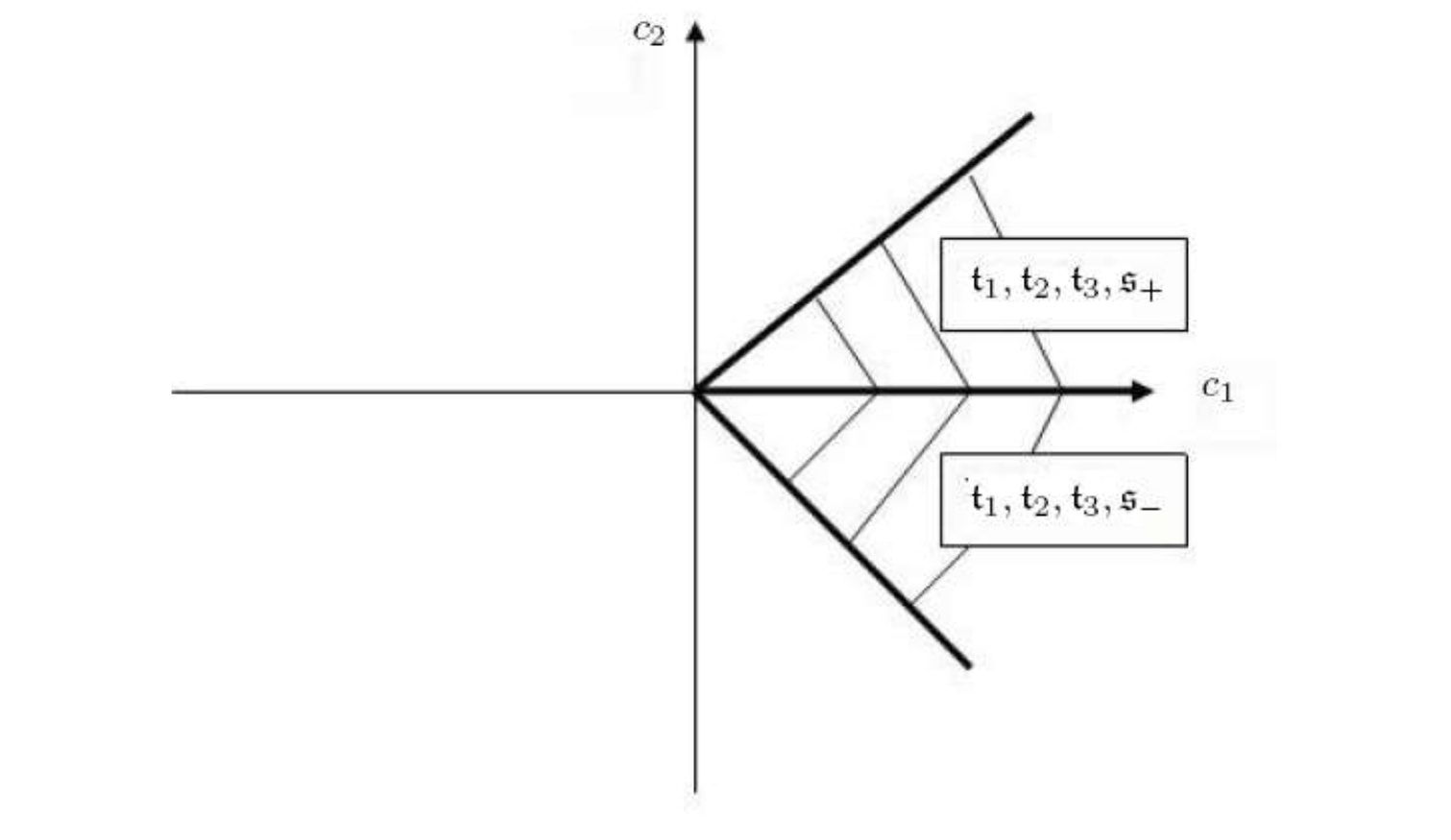}
\end{center}
\caption{\footnotesize
The Weyl chamber. The open wedge corresponds to regular orbits and
the two lines with slopes $\pm \pi/4$ correspond to singular orbits. On a regular orbit corresponding to the strictly upper domain we find equilibria  of type $\mathfrak{t}_1$, $\mathfrak{t}_2$, $\mathfrak{t}_3$, $\mathfrak{s}_{+}$. On a regular orbit corresponding to the strictly lower domain we find equilibria of type $\mathfrak{t}_1$, $\mathfrak{t}_2$, $\mathfrak{t}_3$, $\mathfrak{s}_{-}$. On a regular orbit corresponding to the line $c_2=0$ we find only equilibria of type $\mathfrak{t}_1$, $\mathfrak{t}_2$, $\mathfrak{t}_3$.
}
\label{equilibria}
\end{figure}

\noindent \textbf{Remark.} The free rigid body equations \eqref{soM} on 
$\mathfrak{so}(4)$ should not be confused with another well-known integrable system, also on $\mathfrak{so}(4)$, describing the motion of
a rigid body in an ideal fluid (the Clebsch system). The latter is naturally
a  Hamilton-Poisson system on the Euclidean Lie algebra 
$\mathfrak{se}(3)$;  Bobenko \cite{Bobenko86} (see also \cite{Pogosyan83}, \cite{Pogosyan84}, \cite{Tsiganov06}), realized it also
as a Hamilton-Poisson system on $\mathfrak{so}(4)$.

\section{Nonlinear stability}
In this section we study the nonlinear stability of the
equilibrium states ${\cal E}\cap \operatorname{Orb}_{c_1;c_2}$ for the dynamics \eqref{soM} on a generic adjoint orbit.

Since the system \eqref{soM} on a generic adjoint orbit is completely integrable (\cite{BoFo04}, \cite{Fo88}, \cite{Mishchenko70}, \cite{MiFo78}), for the 
$\mathfrak{so}(4)$ case we have a supplementary constant of motion. The 4th order Mishchenko's constant of motion for the equations \eqref{soM} that commutes with $H$ (see, \cite{Mishchenko70}, \cite{Ratiu80}, \cite{Morosi}) is given by
$$
I(M)=(\lambda_2^2+\lambda_3^2)x_1^2+(\lambda_1^2+\lambda_3^2)x_2^2+(\lambda_1^2+\lambda_2^2)x_3^2+(\lambda_1^2+\lambda_4^2)y_1^2+(\lambda_2^2+\lambda_4^2)y_2^2+(\lambda_3^2+\lambda_4^2)y_3^2.
$$
Without loss of generality, we can choose an ordering for $\lambda_i$'s, namely
$$\lambda_1>\lambda_2>\lambda_3>\lambda_4.$$

The restriction of the dynamics $\dot{M} = [M, \Omega]$ 
to the regular adjoint orbit $\operatorname{Orb}_{c_1;c_2}$ 
is thus a completely integrable Hamiltonian system
\begin{equation}\label{sored}
\left(\operatorname{Orb}_{c_1;c_2},\omega_{\operatorname{Orb}_{c_1;c_2}},H|_{\operatorname{Orb}_{c_1;c_2}}\right),
\end{equation} 
where $\omega_{\operatorname{Orb}_{c_1;c_2}}$ is the orbit
symplectic structure on $\operatorname{Orb}_{c_1;c_2}$.
The Hamiltonian system \eqref{sored} has all equilibria given by Corollary \ref{cor_equ_orbit}. These are of two types:
\begin{align*}
\mathcal{K}_0&: = \left\{M\in \operatorname{Orb}_{c_1;c_2} \mid \mathbf{d}\left(H|_{\operatorname{Orb}_{c_1;c_2}}\right)(M) = 0, \;
\mathbf{d}\left(I|_{\operatorname{Orb}_{c_1;c_2}}\right)(M)=0\right\}\\
\mathcal{K}_1 & : =  \left\{M\in \operatorname{Orb}_{c_1;c_2} \mid \mathbf{d}\left(H|_{\operatorname{Orb}_{c_1;c_2}}\right)(M) = 0, \;
\mathbf{d}\left(I|_{\operatorname{Orb}_{c_1;c_2}}\right)(M) \neq 0\right\}
\end{align*}

By a direct computation, the equilibria of the system \eqref{soM} split in the above two types as follows.
\begin{prop} \label{two_types_equilibria}
$\mathcal{K}_0 = \operatorname{Orb}_{c_1;c_2} \cap (\mathfrak{t}_1\cup\mathfrak{t}_2 \cup \mathfrak{t}_3)$ and $\mathcal{K}_1 = \operatorname{Orb}_{c_1;c_2} \cap \left[( \mathfrak{s}_+ \cup \mathfrak{s}_-)\setminus (\mathfrak{t}_1\cup\mathfrak{t}_2 \cup \mathfrak{t}_3) \right]$.
\end{prop}

\medskip

Since the system \eqref{sored} is completely integrable, we have 
$\{I|_{\operatorname{Orb}_{c_1;c_2}}, H|_{\operatorname{Orb}_{c_1;c_2}}\}=0$, which implies that at an equilibrium $M  \in \operatorname{Orb}_{c_1;c_2}$ we get
$$
\left[\mathbf{D}X_{I|_{\operatorname{Orb}_{c_1;c_2}}}(M), \mathbf{D}X_{H|_{\operatorname{Orb}_{c_1;c_2}}}(M)\right]=0,
$$
where $\mathbf{D}X_{I|_{\operatorname{Orb}_{c_1;c_2}}}(M)$ and $\mathbf{D}X_{H|_{\operatorname{Orb}_{c_1;c_2}}}(M)$ are the derivatives of the vector fields  $X_{I|_{\operatorname{Orb}_{c_1;c_2}}}$ and $X_{H|_{\operatorname{Orb}_{c_1;c_2}}} $ at the equilibrium $M $. Thus $\mathbf{D}X_{I|_{\operatorname{Orb}_{c_1;c_2}}}(M)$, $\mathbf{D}X_{H|_{\operatorname{Orb}_{c_1;c_2}}}(M)$ are infinitesimally symplectic relative to the symplectic form $\omega_{\operatorname{Orb}_{c_1;c_2}}(M)$ on the vector space  $T_M  \operatorname{Orb}_{c_1;c_2}$.

An equilibrium point $M\in \mathcal{K}_0$ is called \emph{non-degenerate} if $\mathbf{D}X_{I|_{\operatorname{Orb}_{c_1;c_2}}}(M)$ and $\mathbf{D}X_{H|_{\operatorname{Orb}_{c_1;c_2}}}(M)$ generate a Cartan subalgebra of the Lie algebra of infinitesimal linear transformations of the symplectic vector space $\left(T_M  \operatorname{Orb}_{c_1;c_2}, \omega_{\operatorname{Orb}_{c_1;c_2}}(M)\right)$

It follows that for a non-degenerate equilibrium belonging to $\mathcal{K}_0$ the matrices $\mathbf{D}X_{I|_{\operatorname{Orb}_{c_1;c_2}}}(M)$ and $\mathbf{D}X_{H|_{\operatorname{Orb}_{c_1;c_2}}}(M)$ can be simultaneously conjugated to one of the following four Cartan subalgebras
\begin{align}
&\text{Type 1:} \quad 
\begin{bmatrix}
0&0&\!\!\!\!-A&0\\
0&0&0&\!\!\!\!-B\\
A&0&0&0\\
0&B&0&0
\end{bmatrix}
\qquad\qquad \text{Type 2:} \quad
\begin{bmatrix}
-A&0&0&0\\
\;\;\;0&0&0&\!\!\!\!-B\\
\;\;\;0&0&A&0\\
\;\;\;0&B&0&0
\end{bmatrix} \nonumber \\ 
\label{cartan_subalgebras}
\\
&\text{Type 3:} \quad 
\begin{bmatrix}
-A&0&0&0\\
\;\;\;0&\!\!\!\!-B&0&0\\
\;\;\;0&0&A&0\\
\;\;\;0&B&0&B
\end{bmatrix}
\qquad\qquad \text{Type 4:} \quad
\begin{bmatrix}
-A&-B&0&\;\;0\\
\;\;\;B&-A&0&\;\;0\\
\;\;\;0&\;\;\;0&A&\!\!-B\\
\;\;\;0&\;\;\;0&B&\;\;A
\end{bmatrix} \nonumber
\end{align}
where $A,B \in \mathbb{R}$ (see, e.g., \cite{BoFo04}, Theorems 1.3 and 1.4).\\
Equilibria of type 1 are called {\it center-center} with the corresponding eigenvalues for the linearized system: $i A,-i A,i B,-i B$.\\
Equilibria of type 2 are called {\it center-saddle} with the corresponding eigenvalues for the linearized system: $A,-A,i B,-iB$.\\
Equilibria of type 3 are called {\it saddle-saddle} with the corresponding eigenvalues for the linearized system: $A,-A,B,-B$.\\
Equilibria of type 4 are called {\it focus-focus} with the corresponding eigenvalues for the linearized system: $A+i B,A-i B,-A+i B,-A-i B$.
\medskip

We will begin the study of stability problem for the three families of equilibria corresponding to the three coordinate-type Cartan subalgebras.\\
\paragraph{The equilibria in $\mathfrak{t}_1$.} We begin with the study of stability and non-degeneracy for the equilibria $M_{a,b}^1\in \mathfrak{t}_1\cap \operatorname{Orb}_{c_1;c_2}$ (see Theorem \ref{Weil_group_orbit}). 

The characteristic polynomial of $\mathbf{D}X_{H|_{\operatorname{Orb}_{c_1;c_2}}}(M_{a,b}^1)$ is
\begin{align*}
z^4 &+ \left[ \frac{a^2}{( \lambda_2+ \lambda_3)^2} \left(\frac{( \lambda_2 - \lambda_4)( \lambda_3- \lambda_4)}{( \lambda_3+ \lambda_4)( \lambda_2+ \lambda_4)} + \frac{( \lambda_1 - \lambda_3)( \lambda_1- \lambda_2)}{( \lambda_1+ \lambda_2)( \lambda_1+ \lambda_3)} \right) \right. \\
& \qquad \qquad  \left.
- \frac{b^2}{( \lambda_1+ \lambda_4)^2} \left(\frac{( \lambda_1 - \lambda_2)( \lambda_2- \lambda_4)}{( \lambda_2+ \lambda_4)( \lambda_1+ \lambda_2)} + 
\frac{( \lambda_1 - \lambda_3)( \lambda_3- \lambda_4)}{( \lambda_3+ \lambda_4)( \lambda_1+ \lambda_3)} \right)
\right]z^2 \\
& + \frac{(\lambda _1- \lambda_2)(\lambda _1- \lambda_3) (\lambda _2- \lambda_4)(\lambda _3- \lambda_4)}{(\lambda _1+ \lambda_2)(\lambda _1+ \lambda_3)(\lambda _3+ \lambda_4)(\lambda _2+ \lambda_4)} \left( \frac{a^2}{(\lambda_2+\lambda_3)^2}-\frac{b^2}{(\lambda_1+\lambda_4)^2} \right)^2.
\end{align*}
The discussion of the position of the four roots in the complex plane is very complicated since the signs of the coefficients of $z^2$ and $z^0$ vary and depend on the relative size of the real numbers $a$ and $b$ which are arbitrary.

Thus, we proceed in a different way. We have already seen that the linear operators $\mathbf{D}X_{H|_{\operatorname{Orb}_{c_1;c_2}}}(M^1_{a,b})$ and
$\mathbf{D}X_{I|_{\operatorname{Orb}_{c_1;c_2}}}(M^1_{a,b})$ are commuting and it is also easy to see that they generate a 2-dimensional subspace.

Therefore, the span of $\mathbf{D}X_{H|_{\operatorname{Orb}_{c_1;c_2}}}(M^1_{a,b})$ and
$\mathbf{D}X_{I|_{\operatorname{Orb}_{c_1;c_2}}}(M^1_{a,b})$ forms a two dimensional Abelian subalgebra of the infinitesimally symplectic linear maps on $\left(T_{M^1_{a,b}} \operatorname{Orb}_{c_1;c_2}, \omega|_{\operatorname{Orb}_{c_1;c_2}} \right)$. We want to show that it is a Cartan subalgebra in order to conclude that that $M^1_{a,b} $ is a non-degenerate equilibrium. This is the case if and only if $\operatorname{span}_ \mathbb{R}\left\{\mathbf{D}X_{H|_{\operatorname{Orb}_{c_1;c_2}}}(M^1_{a,b}),\,
\mathbf{D}X_{I|_{\operatorname{Orb}_{c_1;c_2}}}(M^1_{a,b}) \right\}$ contains an element all of whose eigenvalues are distinct (see, e.g. \cite{BoFo04}, \S1.8.2). 

To show the existence of such an element we begin with the study of the
characteristic polynomial 
$$z^4+v_1z^2+w_1=0,$$
of $\mathbf{D}X_{I|_{\operatorname{Orb}_{c_1;c_2}}}(M^1_{a,b})$,  where
\begin{align*}
w_1&=16(\lambda_1^2-\lambda_2^2)(\lambda_1^2-\lambda_3^2)(\lambda_2^2-\lambda_4^2)(\lambda_3^2-\lambda_4^2)(a^2-b^2)^2>0 \quad \text{since} \quad a \neq b \\
v_1&=S_1a^2+T_1b^2\\
S_1&=4\left(2 \lambda_{2}^{2} \lambda_{3}^{2}-\lambda_{3}^{2}
\lambda_{4}^{2}+ \lambda_{1}^{4}+ \lambda_{4}^{4}- 
\lambda_{2}^2\lambda_{4}^{2}- \lambda_{1}^{2}\lambda_{
2}^{2}- \lambda_{1}^{2} \lambda_{3}^{2}\right) \\
T_1&=4\left(\lambda_{2}^{4}- \lambda_{1}^{2}\lambda_{2}^{2}+2 \lambda_{1}^{2}\lambda_{4}^{2}-\lambda_{2}^{2}\lambda_{4}^{2}-\lambda_{1}^{2}\lambda_{3}^{2}-\lambda_{3}
^{2}\lambda_{4}^{2}+\lambda_{3}^{4}\right).
\end{align*}
Using \eqref{values_a_b} we have 
$$v_1=c_1(S_1+T_1)+\sqrt{c_1^2-c_2^2}(S_1-T_1)$$
and
\begin{align*}
S_1+T_1&=4(\lambda_1^2+\lambda_4^2-\lambda_2^2-\lambda_3^2)^2\geq 0\\
S_1-T_1&=4(\lambda_1^2-\lambda_3^2+\lambda_2^2-\lambda_4^2)(\lambda_1^2-\lambda_2^2+\lambda_3^2-\lambda_4^2)>0
\end{align*}
which shows that $v_1>0$.

The discriminant of the quadratic equation in $z^2$ is
$$
\Delta_1=v_1^2-4w_1=16(\lambda_1^2+\lambda_4^2-\lambda_2^2-\lambda_3^2)^2\left(S'_1a^4+T'_1a^2b^2+U'_1b^4\right),
$$
where
\begin{align*}
S'_1&=(\lambda_1^2-\lambda_4^2)^2>0, \qquad U'_1=(\lambda_2^2-\lambda_3^2)^2>0\\
T'_1&=2\left(-\lambda_3^2 \lambda_4^2 - \lambda_1^2
\lambda_3^2 +2\lambda_1^2 \lambda_4^2 + 2
\lambda_2^2 \lambda_3^2 - \lambda_1^2 \lambda_
2^2 - \lambda_2^2 \lambda_4^2\right).
\end{align*}
Furthermore,
\begin{align*}
S'_1a^4+T'_1a^2b^2+U'_1b^4&=2c_1^2(S'_1+U'_1)+c_2^2(T'_1-S'_1-U'_1)+2c_1\sqrt{c_1^2-c_2^2}(S'_1-U'_1)\\
&=2(c_1^2-c_2^2)(S'_1+U'_1)+c_2^2(S'_1+T'_1+U'_1)+2c_1\sqrt{c_1^2-c_2^2}(S'_1-U'_1).
\end{align*}
Since $\lambda_1 > \lambda_2 > \lambda_3 > \lambda_4$ we have
\begin{align*}
S'_1+U'_1&=(\lambda_1^2-\lambda_4^2)^2+(\lambda_2^2-\lambda_3^2)^2>0\\
S'_1-U'_1&=(\lambda_1^2+\lambda_2^2 -\lambda_3^2 -\lambda_4^2)(\lambda_1^2-\lambda_2^2+\lambda_3^2-\lambda_4^2)>0\\
S'_1+U'_1+T'_1&=(\lambda_1^2+\lambda_4^2-\lambda_2^2-\lambda_3^2)^2\geq 0.
\end{align*}
{\bf Case I.} $\lambda_1^2+\lambda_4^2\not= \lambda_2^2+\lambda_3^2$\\
Thus we have $\Delta_1>0$ (recall $c_1 > |c_2|$). Since $v_1>0,w_1>0$, the equation $t^2+v_1t+w_1=0$ has two non-zero distinct negative real roots and, therefore, the equation $z^4+v_1z^2+w_1=0$ has two distinct pairs of purely imaginary roots different from zero.
Thus  $\operatorname{span}_ \mathbb{R}\left\{\mathbf{D}X_{H|_{\operatorname{Orb}_{c_1;c_2}}}(M^1_{a,b}),\,
\mathbf{D}X_{I|_{\operatorname{Orb}_{c_1;c_2}}}(M^1_{a,b}) \right\}$ is a Cartan subalgebra and it is of the first type in \eqref{cartan_subalgebras}. It follows that the equilibrium $M_{a,b}^1$ is non-degenerate and nonlinearly stable because it is of center-center type (see \cite{BoFo04}, Theorem 1.5). 

The above computations are identical if one replaces $a$ by
$-a$ and $b$ by $-b$. So the same argument applies to 
$M^1_{-a, -b}$.

\medskip

\noindent {\bf Case II.} $\lambda_1^2+\lambda_4^2= \lambda_2^2+\lambda_3^2$\\
The eigenvalues of $\mathbf{D}X_{I|_{\operatorname{Orb}_{c_1;c_2}}}(M^1_{a,b})$ are conjugate purely imaginary of multiplicity two. In order to determine non-degeneracy of the equilibrium $M_{a,b}^1$ we have to find a linear combination $\mathbf{D}X_{H|_{\operatorname{Orb}_{c_1;c_2}}}(M^1_{a,b}) +\alpha \mathbf{D}X_{I|_{\operatorname{Orb}_{c_1;c_2}}}(M^1_{a,b})$, where $\alpha$ is a non-zero real number, that has distinct eigenvalues. The eigenvalues of this linear combination are the roots of the equation
$$u_1'z^4+v_1'z^2+w_1'=0,$$
where
\begin{align*}
u_1'&=(\lambda_1+\lambda_2)(\lambda_1+\lambda_3)(\lambda_1+\lambda_4)^4(\lambda_2+\lambda_3)^4(\lambda_2+\lambda_4)(\lambda_3+\lambda_4)>0\\
w_1'&=(\lambda_1-\lambda_2)(\lambda_1-\lambda_3)(\lambda_2-\lambda_4)(\lambda_3-\lambda_4)(X_1\alpha^2+Y_1\alpha+Z_1)^2 \geq 0,
\end{align*}
with
$$X_1=4(\lambda_1+\lambda_2)(\lambda_1+\lambda_3)(\lambda_1+\lambda_4)^2(\lambda_2+\lambda_3)^2(\lambda_2+\lambda_4)(\lambda_3+\lambda_4)(a^2-b^2)\not= 0$$
and $Y_1,Z_1$ are also expressions of $\lambda_1,\lambda_2,\lambda_3,\lambda_4,a,b$.

The discriminant of the quadratic equation $u_1't^2+v_1't+w_1'=0$ obtained by denoting $z^2=t$ is 
\[
\Delta_1'=4(\lambda_1+\lambda_4)^6(\lambda_2+\lambda_3)^6(Y_2\alpha+Z_2)^2D,
\] 
where 
\begin{align*}
D&=2\left[(c_1^2-c_2^2)\left((\lambda_1^2-\lambda_4^2)^2+(\lambda_2^2-\lambda_3^2)^2\right)+c_1\sqrt{c_1^2-c_2^2}(\lambda_1^2-\lambda_3^2+\lambda_2^2-\lambda_4^2)(\lambda_1^2-\lambda_2^2+\lambda_3^2-\lambda_4^2)\right]>0 \\
Y_2&=-2(\lambda_1+\lambda_2)(\lambda_1+\lambda_3)(\lambda_2+\lambda_4)(\lambda_3+\lambda_4)(\lambda_1+\lambda_4-\lambda_2-\lambda_3) \\
Z_2&=\lambda_1\lambda_4-\lambda_2\lambda_3.
\end{align*}
Note that $Y_2 = 0 $ if and only if $\lambda_1+\lambda_4= \lambda_2+
\lambda_3$, because $\lambda_i + \lambda_j >0$ for any $i \neq j $. Then 
$(\lambda_1+\lambda_4)^2= (\lambda_2+\lambda_3)^2$ and  $
\lambda_1^2+\lambda_4^2= \lambda_2^2+\lambda_3^2$ imply $Z_2 = 0 $.  Conversely, if $Z_2 = 0 $, then $\lambda_1^2+\lambda_4^2= 
\lambda_2^2+\lambda_3^2$ implies that $(\lambda_1+\lambda_4)^2= 
(\lambda_2+\lambda_3)^2$, that is $\lambda_1+\lambda_4= \pm 
(\lambda_2+\lambda_3)$. Since the solution with minus is not possible 
because $\lambda_i + \lambda_j >0$ for any $i \neq j $, we conclude that 
$Y_2 = 0 $. Thus, $Y_2 = 0 $ if and only if $Z_2 = 0 $. 

However, if $Y_2 = 0$, so $Z_2=0$ which means that $\lambda_1 \lambda_2 = \lambda_3 \lambda_4 $, then we also have $(\lambda_1 - \lambda_4)^2 
=(\lambda_2 - \lambda_3)^2$ and consequently 
$\lambda_1 - \lambda_4 =\pm(\lambda_2 - \lambda_3)$. The solution with 
minus is impossible because $\lambda _1 + \lambda_2 > \lambda_3 + 
\lambda_4 $ since, by hypothesis, $\lambda_1 > \lambda_2 > 
\lambda_3 > \lambda_4$. Hence we must have
$\lambda_1 - \lambda_4=\lambda_2 - \lambda_3$ which together with $
\lambda_1+\lambda_4= \lambda_2+ \lambda_3$ implies that 
$\lambda_1=\lambda_2$ which is also impossible since $\lambda_1> \lambda_2$. Therefore  $Y_2 \neq 0$ and hence $\Delta_1' >0$  if we choose $\alpha\neq -Z_2/Y_2$.

Furthermore, $v_1'$ has the expression
$$v_1'=2(\lambda_1+\lambda_4)^2(\lambda_2+\lambda_3)^2\sqrt{c_1^2-c_2^2}(X_3\alpha^2+Y_3\alpha+Z_3),$$
where
$$X_3=2(\lambda_1+\lambda_2)(\lambda_1+\lambda_3)(\lambda_1+\lambda_4)^2(\lambda_2+\lambda_3)^2(\lambda_2+\lambda_4)(\lambda_3+\lambda_4)
(\lambda_1^2-\lambda_3^2+\lambda_2^2-\lambda_4^2)(\lambda_1^2-\lambda_2^2+\lambda_3^2-\lambda_4^2)>0.$$

Since $u_1' >0 $ and $\Delta_1' >0 $ for $\alpha\neq -Z_2/Y_2$, choosing $\alpha \in \mathbb{R} $ large enough we also have $v_1'>0$, $w_1'>0$, and hence the matrix  
$\mathbf{D}X_{H_{|_{\operatorname{Orb}_{c_1;c_2}}}}(M_{a,b}^1)+\alpha \mathbf{D}X_{I_{|_{\operatorname{Orb}_{c_1;c_2}}}}(M_{a,b}^1)$ has four distinct purely imaginary eigenvalues (in particular, zero is not an eigenvalue). Therefore $M_{a,b}^1$ is 
a non-degenerate equilibrium and $\operatorname{span}_ \mathbb{R}\left\{\mathbf{D}X_{H|_{\operatorname{Orb}_{c_1;c_2}}}(M^1_{a,b}),\,
\mathbf{D}X_{I|_{\operatorname{Orb}_{c_1;c_2}}}(M^1_{a,b}) \right\}$ is a Cartan subalgebra of the first type in 
\eqref{cartan_subalgebras}. The equilibrium $M_{a,b}^1$ is 
thus of center-center type and therefore it is nonlinearly 
stable on the adjoint orbit of $\mathfrak{so}(4)$ determined
by $a$ and $b$.

Since $M_{a,b}^1$ lies on a generic adjoint orbit, the
analysis above showing stability on the adjoint orbit determined by $a,b$, implies that this equilibrium is also
nonlinearly stable for the Lie-Poisson dynamics on $\mathfrak{so}(4)$. The same holds for the equilibrium $M_{-a,-b}^1$. 

\begin{thm}
The equilibria $M_{a,b}^1, M_{-a,-b}^1 \in 
\mathfrak{t}_1 \cap \operatorname{Orb}_{c_1; c_2}$ are non-degenerate of 
center-center type and therefore nonlinearly stable on the corresponding adjoint orbit and also nonlinearly stable
for the Lie-Poisson dynamics on $\mathfrak{so}(4)$.
\end{thm}

Next, we study the stability of the equilibria $M_{b,a}^1, 
M_{-b,-a}^1 \in \mathfrak{t}_1$. The characteristic 
polynomial of $\mathbf{D}X_{I|\operatorname{Orb}_{c_1;c_2}}
(M^1_{b,a})$ relative to the matrix $[\omega(M^1_{b,a})]$ is
$z^4+\widetilde{v}_1z^2+\widetilde{w}_1=0$, where 
\[
\widetilde{v}_1=\widetilde{S}_1a^2+\widetilde{T}_1b^2
\]
with
$$\widetilde{S}_1=-4[(\lambda_1^2-\lambda_2^2)(\lambda_2^2-\lambda_4^2)+(\lambda_1^2-\lambda_3^2)(\lambda_3^2-\lambda_4^2)]<0;$$
$$\widetilde{T}_1=4[(\lambda_1^2-\lambda_2^2)(\lambda_1^2-\lambda_3^2)+(\lambda_2^2-\lambda_4^2)(\lambda_3^2-\lambda_4^2)]>0$$
and
$$\widetilde{w}_1=16(\lambda_1^2-\lambda_2^2)(\lambda_1^2-\lambda_3^2)(\lambda_2^2-\lambda_4^2)(\lambda_3^2-\lambda_4^2)(a^2-b^2)^2>0.$$
The discriminant of the quadratic equation associated to the characteristic equation is:
$$
\widetilde{\Delta}_1=\widetilde{v}_1^2-4\widetilde{w}_1
= 16(\lambda_1^2- \lambda_2^2- \lambda_3^2+ \lambda_4^2)^2
a^4\left[\widetilde{S}_2\left(\frac{b^2}{a^2}\right)^2+\widetilde{T}_2\frac{b^2}{a^2}+\widetilde{U}_2\right],$$
where
$$\widetilde{S}_2=(\lambda_1^2-\lambda_4^2)^2>0;~~~\widetilde{U}_2=(\lambda_2^2-\lambda_3^2)^2>0;$$
$$\widetilde{T}_2=-2[(\lambda_1^2-\lambda_2^2)(\lambda_3^2-\lambda_4^2)+(\lambda_1^2-\lambda_3^2)(\lambda_2^2-\lambda_4^2)]<0.$$
We introduce the following quadratic function related to $\widetilde{\Delta}_1$:
$$\widetilde{f}(t)=\widetilde{S}_2t^2+\widetilde{T}_2t+\widetilde{U}_2.$$
The discriminant of $\widetilde{f}$ is given by
$$\Delta_{\widetilde{f}}=16(\lambda_1^2-\lambda_2^2)(\lambda_1^2-\lambda_3^2)(\lambda_2^2-\lambda_4^2)(\lambda_3^2-\lambda_4^2)>0.$$
Since $\Delta_{\widetilde{f}}>0,-\frac{\widetilde{T}_2}{\widetilde{S}_2}>0,\frac{\widetilde{U}_2}{\widetilde{S}_2}>0$, the quadratic equation associated to $\widetilde{f}$ has two distinct strictly positive real solutions, which we will denote by $\alpha_1<\alpha_2$.

Further we will distinguish two cases, when $\lambda_1^2+ \lambda_4^2\not= \lambda_2^2+\lambda_3^2$ and $\lambda_1^2+ \lambda_4^2= \lambda_2^2+\lambda_3^2$.\\\\
{\bf Case I.} $\lambda_1^2+ \lambda_4^2\not= \lambda_2^2+\lambda_3^2$\\
We notice that $\widetilde{f}(1)=(\lambda_1^2- \lambda_2^2- \lambda_3^2+ \lambda_4^2)^2>0$ and
$$-\frac{\widetilde{T}_2}{2\widetilde{S}_2}-1=-\frac{(\lambda_1^2-\lambda_2^2)(\lambda_1^2-\lambda_3^2)+(\lambda_2^2-\lambda_4^2)(\lambda_3^2-\lambda_4^2)}{(\lambda_1^2-\lambda_4^2)^2}<0,$$ which implies the following ordering
$$0<\alpha_1<\alpha_2<1.$$
We denote by $\alpha_3:=-\frac{\widetilde{S}_1}{\widetilde{T}_1}$.
Also
$$\widetilde{f}(\alpha_3)=-\frac{64(\lambda_1^2- \lambda_2^2- \lambda_3^2+ \lambda_4^2)^2(\lambda_1^2-\lambda_2^2)(\lambda_1^2-\lambda_3^2)(\lambda_2^2-\lambda_4^2)(\lambda_3^2-\lambda_4^2)}{\widetilde{T}_1^2}<0.$$
Consequently, we further have the following ordering
$$0<\alpha_1<\alpha_3<\alpha_2<1.$$
We introduced $\alpha_3$ as the value where $\widetilde{v}_1$ changes sign. More precisely, we notice that $\widetilde{v}_1
\gtrless 0$ if and only if $\frac{b^2}{a^2} \gtrless \alpha_3$.
We will distinguish the following subcases:\\\\
{\bf Subcase 1.} $\frac{b^2}{a^2}\in [0,\alpha_1)$\\
In this situation we have $\widetilde{\Delta}_1>0$ and $\widetilde{v}_1<0$. Consequently, the eigenvalues of $\mathbf{D}X_{I|\operatorname{Orb}_{c_1;c_2}}
(M^1_{b,a})$ are of the form $A,-A,B,-B$, with $A\not= B; A, B\not= 0$, and $M^1_{b,a}$ is {\it unstable of saddle-saddle type}.\\
{\bf Subcase 2.} $\frac{b^2}{a^2}\in (\alpha_1,\alpha_3)\cup (\alpha_3,\alpha_2)$\\
In this situation we have $\widetilde{\Delta}_1<0$ and $\widetilde{v}_1\not= 0$. Consequently, the eigenvalues of $\mathbf{D}X_{I|\operatorname{Orb}_{c_1;c_2}}
(M^1_{b,a})$ are of the form $A+iB,A-iB,-A+iB,-A-iB$, with $A\not= B; A, B\not= 0$, and $M^1_{b,a}$ is {\it unstable of focus-focus type}.\\
{\bf Subcase 3.} $\frac{b^2}{a^2}\in (\alpha_2,1)$\\
In this situation we have $\widetilde{\Delta}_1>0$ and $\widetilde{v}_1>0$. Consequently, the eigenvalues of $\mathbf{D}X_{I|\operatorname{Orb}_{c_1;c_2}}
(M^1_{b,a})$ are of the form $iA,-iA,iB,-iB$, with $A\not= B; A, B\not= 0$, and $M^1_{b,a}$ is {\it stable of center-center type}.\\
\noindent {\bf Subcase 4.} $\frac{b^2}{a^2}= \alpha_3$\\
In this situation we have $\widetilde{\Delta}_1>0$ and $\widetilde{v}_1=0$. Consequently, the eigenvalues of $\mathbf{D}X_{I|\operatorname{Orb}_{c_1;c_2}}
(M^1_{b,a})$ are of the form $A+iA,A-iA,-A+iA,-A-iA$, with $A\not= 0$, and $M^1_{b,a}$ is {\it unstable of focus-focus type}.

\bigskip

\noindent The following subcases are frontier cases.
\bigskip

\noindent{\bf Subcase 5.} $\frac{b^2}{a^2}= \alpha_1$\\
In this situation we have $\widetilde{\Delta}_1=0$ and $\widetilde{v}_1<0$. Consequently, the eigenvalues of $\mathbf{D}X_{I|\operatorname{Orb}_{c_1;c_2}}
(M^1_{b,a})$ are of the form $A,-A,A,-A$, with $A\not= 0$, and $M^1_{b,a}$ is {\it unstable}.\\
To decide the type of instability, we need to determine the non-degeneracy of $M^1_{b,a}$, i.e. we have to find a linear combination $\mathbf{D}X_{H|_{\operatorname{Orb}_{c_1;c_2}}}(M^1_{a,b}) +\beta \mathbf{D}X_{I|_{\operatorname{Orb}_{c_1;c_2}}}(M^1_{b,a})$, where $\beta$ is a non-zero real number, that has distinct eigenvalues. The eigenvalues of this linear combination are the roots of the equation
$$pz^4+qz^2+r=0.$$
The discriminant of the quadratic function associated with the above equation is of the form
$$\Delta=16(X\beta^2+Y\beta+Z)^2a^4\widetilde{f}\left(\frac{b^2}{a^2}\right).$$
We notice that $\Delta =0$ as we are in the case when $\frac{b^2}{a^2}=\alpha_1$ and $\alpha_1$ is a square root of $\widetilde{f}$. Consequently, the linear combination does not have four distinct eigenvalues and so $M^1_{b,a}$ in this frontier subcase is an {\it unstable degenerate equilibrium}.\\
{\bf Subcase 6.} $\frac{b^2}{a^2}= \alpha_2$\\
In this situation we have $\widetilde{\Delta}_1=0$ and $\widetilde{v}_1>0$. Consequently, the eigenvalues of $\mathbf{D}X_{I|\operatorname{Orb}_{c_1;c_2}}
(M^1_{b,a})$ are of the form $iA,-iA,iA,-iA$, with $A\not= 0$. \\
To decide the type of stability, we proceed as in the subcase 5; we need to determine the non-degeneracy of $M^1_{b,a}$, i.e. we have to find a linear combination $\mathbf{D}X_{H|_{\operatorname{Orb}_{c_1;c_2}}}(M^1_{a,b}) +\beta \mathbf{D}X_{I|_{\operatorname{Orb}_{c_1;c_2}}}(M^1_{b,a})$, where $\beta$ is a non-zero real number, that has distinct eigenvalues. The eigenvalues of this linear combination are the roots of the equation
$$pz^4+qz^2+r=0.$$
The discriminant of the quadratic function associated with the above equation is of the form
$$\Delta=16(X\beta^2+Y\beta+Z)^2a^4\widetilde{f}\left(\frac{b^2}{a^2}\right).$$
We notice that $\Delta =0$ as we are in the case when $\frac{b^2}{a^2}=\alpha_2$ and $\alpha_2$ is a square root of $\widetilde{f}$. Consequently, the linear combination does not have four distinct eigenvalues and so $M^1_{b,a}$ in this frontier subcase is a degenerate equilibrium for which we cannot solve the stability problem on the adjoint orbit of
$\mathfrak{so}(4)$ determined by $a$ and $b$.

\begin{thm}
Under the hypothesis $\lambda_1^2+ \lambda_4^2\not= \lambda_2^2+\lambda_3^2$ the following hold:
\begin{itemize}
\item[{\rm (i)}] If $\frac{b^2}{a^2}\in [0,\alpha_1)$, then the equilibria $M_{b,a}^1, M_{-b,-a}^1 \in 
\mathfrak{t}_1 \cap \operatorname{Orb}_{c_1; c_2}$ are non-degenerate unstable of saddle-saddle type on the adjoint orbit determined by $a $ and $b $ and hence are also unstable for the Lie-Poisson dynamics on $\mathfrak{so}(4)$.
\item[{\rm (ii)}] If $\frac{b^2}{a^2}\in (\alpha_1,\alpha_2)$, then the equilibria $M_{b,a}^1, M_{-b,-a}^1 \in 
\mathfrak{t}_1 \cap \operatorname{Orb}_{c_1; c_2}$ are non-degenerate unstable of focus-focus type on the adjoint orbit determined by $a $ and $b $ and hence are also unstable for the Lie-Poisson dynamics on $\mathfrak{so}(4)$.
\item[{\rm (iii)}] If $\frac{b^2}{a^2}\in (\alpha_2,1)$, then the equilibria $M_{b,a}^1, M_{-b,-a}^1 \in 
\mathfrak{t}_1 \cap \operatorname{Orb}_{c_1; c_2}$ are non-degenerate stable of center-center type on the adjoint orbit determined by $a $ and $b $ and hence are also nonlinearly stable for the Lie-Poisson dynamics on $\mathfrak{so}(4)$.
\item[{\rm (iv)}] If $\frac{b^2}{a^2}=\alpha_1$, then the equilibria $M_{b,a}^1, M_{-b,-a}^1 \in 
\mathfrak{t}_1 \cap \operatorname{Orb}_{c_1; c_2}$ are degenerate and unstable on the adjoint orbit determined by $a $ and $b $ and hence are also unstable for the Lie-Poisson dynamics on $\mathfrak{so}(4)$.
\item[{\rm (v)}] If $\frac{b^2}{a^2}=\alpha_2$, then the equilibria $M_{b,a}^1, M_{-b,-a}^1 \in 
\mathfrak{t}_1 \cap \operatorname{Orb}_{c_1; c_2}$ are degenerate and the stability problem on the adjoint orbit determined by $a $ and $b $ remains open. However, these
equilibria are unstable for the Lie-Poisson dynamics on $\mathfrak{so}(4)$.
\end{itemize}
\end{thm}

\begin{figure}[!th]
\begin{center}
\includegraphics[scale=0.6,angle=0]{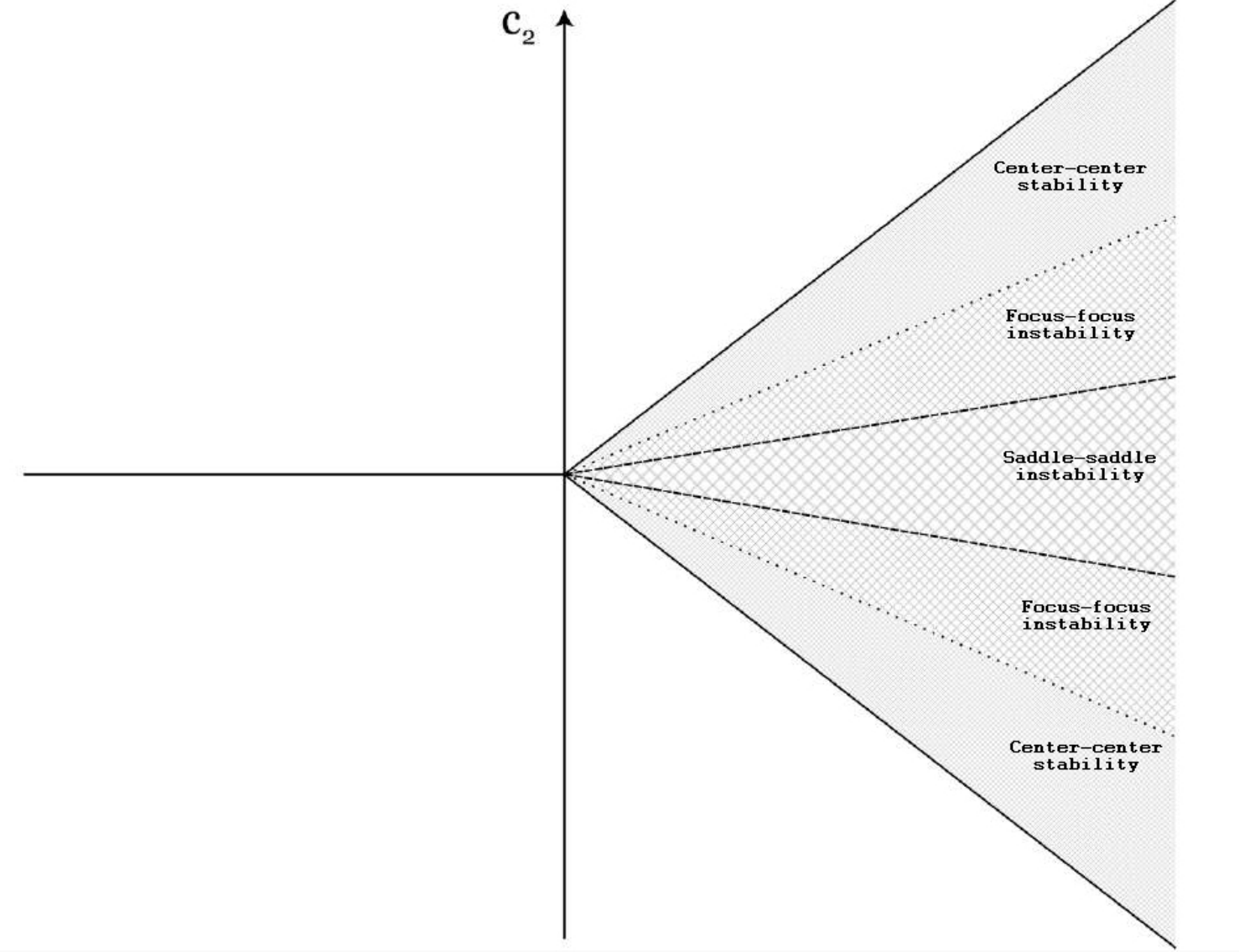}
\end{center}
\caption{\footnotesize
The domain filled with the tight grid corresponds to non-degenerate stable equilibria of center-center type. The domain filled with the medium grid corresponds to non-degenerate unstable equilibria of focus-focus type. The middle wedge filled with the wide grid corresponds to non-degenerate unstable equilibria of saddle-saddle type.
}
\label{figura2}
\end{figure}

The equality $\frac{b^2}{a^2}=\alpha_1$ is equivalent to
$|c_2|=\frac{2\sqrt{\alpha_1}}{1+\alpha_1}c_1$ and thus the coefficients $\pm \frac{2\sqrt{\alpha_1}}{1+\alpha_1}$ represent the slopes of the frontier lines corresponding to case {\rm (iv)} and which are represented by dashed lines in the figure below. The equality $\frac{b^2}{a^2}=\alpha_2$ is equivalent to
$|c_2|=\frac{2\sqrt{\alpha_2}}{1+\alpha_2}c_1$ and thus the coefficients $\pm \frac{2\sqrt{\alpha_2}}{1+\alpha_2}$ represent the slopes of the frontier lines corresponding to case {\rm (v)} and which are represented by dotted lines in the Figure \ref{figura2}.  

We comment on points (iii) and (v).  The statement about the nonlinear
stability in (iii) follows since these equilibria are of center-center type and lie in
the open set of generic adjoint orbits whose intersection with the positive Weyl chamber is the tight grid domain in Figure \ref{figura2}. The instability statement in (v) is due to
the fact that any neighborhood of this equilibrium contains
unstable equilibria of focus-focus type as described in (ii).

Note that case (v) is also the transition point in a Hamiltonian
Hopf bifurcation when the dynamics changes from a center-center configuration (case (iii)) to a focus-focus configuration (case (ii)).
\newpage

\medskip

\noindent {\bf Case II.} $\lambda_1^2+ \lambda_4^2= \lambda_2^2+\lambda_3^2$\\
From the above computations we immediately obtain for this case that $\widetilde{f}(1)=\widetilde{f}(\alpha_3)=0$. Additionally, we notice that
$$1-\alpha_3=4\frac { \left( {\lambda_{{1}}}^{2}-{\lambda_{{2}}}^{2}-{\lambda_{{3}}
}^{2}+{\lambda_{{4}}}^{2} \right) ^{2}}{\widetilde{T}_1}=0.$$
Thus we have the ordering $0<\alpha_1<\alpha_2=\alpha_3=1.$\\
In this case $\widetilde{\Delta}_1=0$ and since $\frac{b^2}{a^2}<\alpha_3=1$ we have $\widetilde{v}_1<0$. Consequently, the eigenvalues of $\mathbf{D}X_{I|\operatorname{Orb}_{c_1;c_2}}
(M^1_{b,a})$ are of the form $A,-A,A,-A$, with $A\not= 0$, and $M^1_{b,a}$ is {\it unstable}.\\

In order to determine non-degeneracy of the equilibrium $M_{b,a}^1$ in this case we have to find a linear combination $\mathbf{D}X_{H|_{\operatorname{Orb}_{c_1;c_2}}}(M^1_{b,a}) +\beta \mathbf{D}X_{I|_{\operatorname{Orb}_{c_1;c_2}}}(M^1_{b,a})$, where $\beta$ is a non-zero real number, that has distinct eigenvalues. The eigenvalues of this linear combination are the roots of the equation
$$\widetilde{u}_1'z^4+\widetilde{v}_1'z^2+\widetilde{w}_1'=0,$$
where
\begin{align*}
\widetilde{u}_1'&=(\lambda_1+\lambda_2)(\lambda_1+\lambda_3)(\lambda_1+\lambda_4)^4(\lambda_2+\lambda_3)^4(\lambda_2+\lambda_4)(\lambda_3+\lambda_4)>0\\
\widetilde{w}_1'&=(\lambda_1-\lambda_2)(\lambda_1-\lambda_3)(\lambda_2-\lambda_4)(\lambda_3-\lambda_4)(Q_1\beta^2+Q_2\beta+Q_3)^2 \geq 0,
\end{align*}
with
$$Q_1=4(\lambda_1+\lambda_2)(\lambda_1+\lambda_3)(\lambda_1+\lambda_4)^2(\lambda_2+\lambda_3)^2(\lambda_2+\lambda_4)(\lambda_3+\lambda_4)(a^2-b^2)\not= 0$$
and $Q_2,Q_3$ are also expressions of $\lambda_1,\lambda_2,\lambda_3,\lambda_4,a,b$.

The discriminant of the quadratic equation $\widetilde{u}_1't^2+\widetilde{v}_1't+\widetilde{w}_1'=0$ obtained by denoting $z^2=t$ is 
\[
\widetilde{\Delta}_1'=4(\lambda_1+\lambda_4)^6(\lambda_2+\lambda_3)^6(P_2\beta+P_3)^2a^4\widetilde{f}\left(\frac{b^2}{a^2}\right),
\] 
where 
\begin{align*}
P_2&=-2(\lambda_1+\lambda_2)(\lambda_1+\lambda_3)(\lambda_2+\lambda_4)(\lambda_3+\lambda_4)(\lambda_1+\lambda_4-\lambda_2-\lambda_3) \\
P_3&=\lambda_1\lambda_4-\lambda_2\lambda_3.
\end{align*}
In this case, by an analogous reasoning as for the equilibria $M^1_{a,b}$, we cannot have simultaneously $P_2=P_3=0$. Consequently, the expression $(P_2\beta+P_3)^2$ in the above discriminant is strictly positive for $\beta$ large enough.\\
Furthermore, $\widetilde{v}_1'$ has the expression
$$\widetilde{v}_1'=-2(\lambda_1+\lambda_4)^2(\lambda_2+\lambda_3)^2(R_1\beta^2+R_2\beta+R_3),$$
where
$$R_1=\frac{1}{4}(\lambda_1+\lambda_2)(\lambda_1+\lambda_3)(\lambda_1+\lambda_4)^2(\lambda_2+\lambda_3)^2(\lambda_2+\lambda_4)(\lambda_3+\lambda_4)a^2\widetilde{T}_1\left(1-\frac{b^2}{a^2}\right)$$
and $R_2,R_3$ are also expressions of $\lambda_1,\lambda_2,\lambda_3,\lambda_4,a,b$.
Since $\frac{b^2}{a^2}<1$, for $\beta$ large enough we obtain $\widetilde{v}_1'<0$.\\
In this case, we further distinguish three subcases:\\\\
{\bf Subcase 1$^{\prime}$.} $\frac{b^2}{a^2}\in [0,\alpha_1)$\\
In this situation we have $\widetilde{\Delta}_1'>0$, $\widetilde{v}_1'<0$ and $\widetilde{w}_1'>0$ for $\beta$ large enough. Consequently, the eigenvalues of $\mathbf{D}X_{I|\operatorname{Orb}_{c_1;c_2}}
(M^1_{b,a})$ are of the form $A,-A,B,-B$, with $A\not= B; A, B\not= 0$, and $M^1_{b,a}$ is {\it unstable of type saddle-saddle}.\\
{\bf Subcase 2$^{\prime}$.} $\frac{b^2}{a^2}\in (\alpha_1,1)$\\ 
In this situation we have $\widetilde{\Delta}_1'<0$, $\widetilde{v}_1'<0$ and $\widetilde{w}_1'>0$ for $\beta$ large enough. Consequently, the eigenvalues of $\mathbf{D}X_{I|\operatorname{Orb}_{c_1;c_2}}
(M^1_{b,a})$ are of the form $A+iB,A-iB,-A+iB,-A-iB$, with $A\not= B; A, B\not= 0$, and $M^1_{b,a}$ is {\it unstable of type focus-focus}.

\bigskip

\noindent We have the following frontier case.
\bigskip

{\bf Subcase 3$^{\prime}$.} $\frac{b^2}{a^2}= \alpha_1$\\
In this situation we have $\widetilde{\Delta}_1'=0$, $\widetilde{v}_1'<0$ and $\widetilde{w}_1'>0$ for $\beta$ large enough. Consequently, the linear combination does not have four distinct eigenvalues and so $M^1_{b,a}$ in this frontier subcase is an {\it unstable degenerate equilibrium}.\\
The same analysis holds for the equilibrium $M^1_{-b,-a}$.

\begin{thm}
Under the hypothesis $\lambda_1^2+ \lambda_4^2= \lambda_2^2+\lambda_3^2$ the following holds:
\begin{itemize}
\item[{\rm (i)}] If $\frac{b^2}{a^2}\in [0,\alpha_1)$, then the equilibria $M_{b,a}^1, M_{-b,-a}^1 \in 
\mathfrak{t}_1 \cap \operatorname{Orb}_{c_1; c_2}$ are non-degenerate unstable of saddle-saddle type.
\item[{\rm (ii)}] If $\frac{b^2}{a^2}\in (\alpha_1,1)$, then the equilibria $M_{b,a}^1, M_{-b,-a}^1 \in 
\mathfrak{t}_1 \cap \operatorname{Orb}_{c_1; c_2}$ are non-degenerate unstable of focus-focus type.  
\item[{\rm (iii)}] If $\frac{b^2}{a^2}=\alpha_1$, then the equilibria $M_{b,a}^1, M_{-b,-a}^1 \in 
\mathfrak{t}_1 \cap \operatorname{Orb}_{c_1; c_2}$ are degenerate and unstable.
\end{itemize}
Thus these equilibria are also unstable for the Lie-Poisson dynamics on $\mathfrak{so}(4)$.
\end{thm}

\begin{figure}[!th]
\begin{center}
\includegraphics[scale=0.6,angle=0]{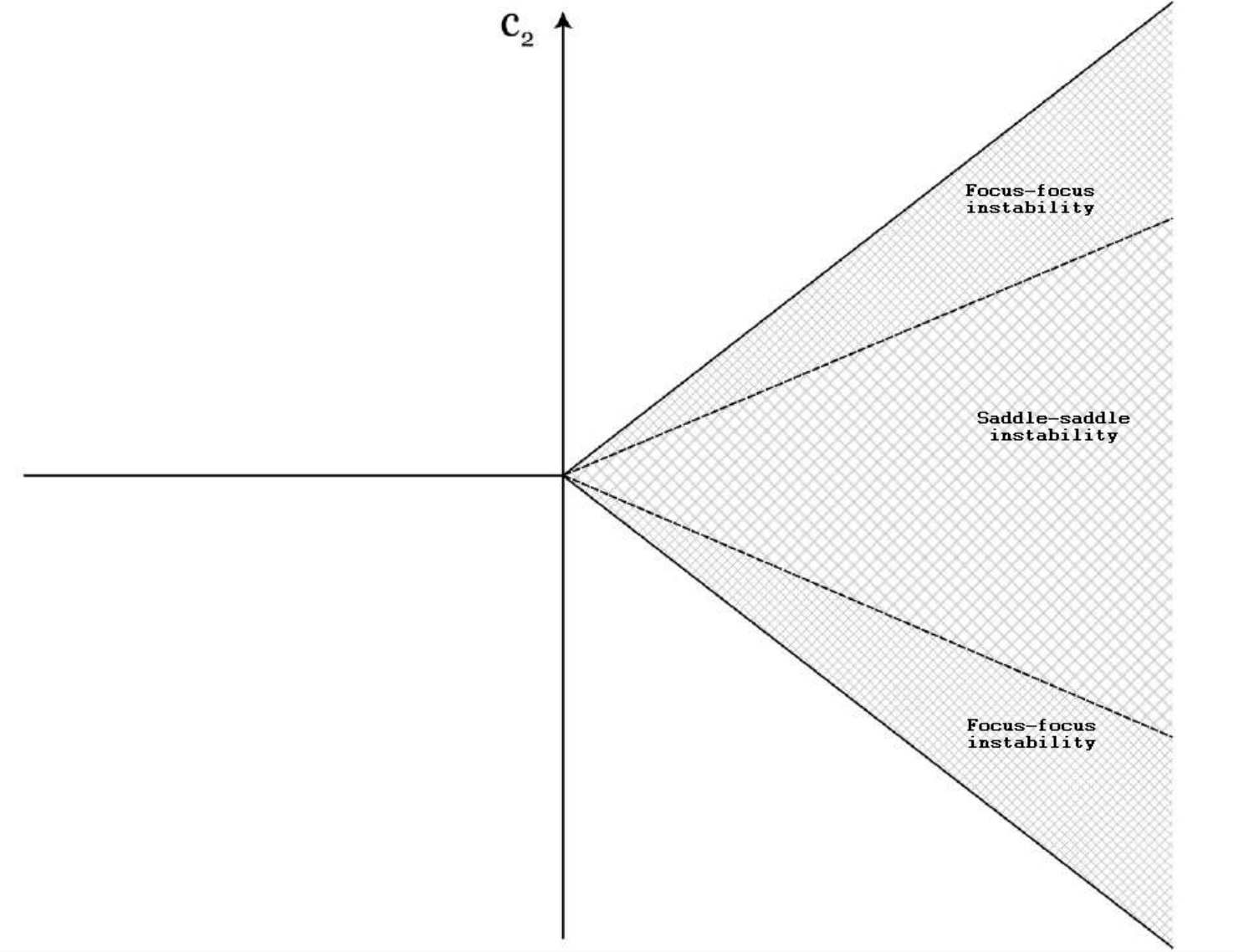}
\end{center}
\caption{\footnotesize
The domain filled with the medium grid corresponds to non-degenerate unstable equilibria of focus-focus type. The domain filled with the wide grid corresponds to non-degenerate unstable equilibria of saddle-saddle type.
}
\label{figura_3}
\end{figure}

The equality $\frac{b^2}{a^2}=\alpha_1$ is equivalent with
$|c_2|=\frac{2\sqrt{\alpha_1}}{1+\alpha_1}c_1$, and thus the coefficients $\pm \frac{2\sqrt{\alpha_1}}{1+\alpha_1}$ represent the slopes of the frontier lines corresponding to case {\rm (iii)} and which are represented by dashed lines in the Figure \ref{figura_3}.

\medskip

\paragraph{The equilibria in $\mathfrak{t}_2$.} We proceed using the same techniques as in the previous case. It is easy to see that $\mathbf{D}X_{I|_{\operatorname{Orb}_{c_1;c_2}}}(M_{a,b}^2)$ and $\mathbf{D}X_{H|_{\operatorname{Orb}_{c_1;c_2}}}(M_{a,b}^2)$ generate a 2-dimensional subspace. The eigenvalues of $\mathbf{D}X_{I|_{\operatorname{Orb}_{c_1;c_2}}}(M_{a,b}^2)$ are the roots of the equation
$$z^4+v_2z^2+w_2=0,$$
where
\begin{align*}
w_2&=-16(\lambda_1^2-\lambda_2^2)(\lambda_1^2-\lambda_4^2)(\lambda_2^2-\lambda_3^2)(\lambda_3^2-\lambda_4^2)(a^2-b^2)^2<0 \\
v_2&=4\left( -\,{\lambda_{{3}}}^{2}{\lambda_{{4}}}^{2}-\,{\lambda_{{2}}}^{2}{\lambda_{{3}}}^{2}-\,{\lambda_{{1}}}^{2}{\lambda_{{4}}}^{2}+\,{
\lambda_{{2}}}^{4}-\,{\lambda_{{1}}}^{2}{\lambda_{{2}}}^{2}+\,{
\lambda_{{4}}}^{4}+2\,{\lambda_{{1}}}^{2}{\lambda_{{3}}}^{2} \right) {
a}^{2}\\
&\quad +4\left( -\,{\lambda_{{3}}}^{2}{\lambda_{{4}}}^{2}+2\,{\lambda_
{{2}}}^{2}{\lambda_{{4}}}^{2}-\,{\lambda_{{1}}}^{2}{\lambda_{{2}}}^{2
}-\,{\lambda_{{1}}}^{2}{\lambda_{{4}}}^{2}+\,{\lambda_{{1}}}^{4}+\,
{\lambda_{{3}}}^{4}-\,{\lambda_{{2}}}^{2}{\lambda_{{3}}}^{2} \right) 
{b}^{2}.
\end{align*}
The quadratic equation $t^2+v_2t+w_2=0$ has the discriminant
$$\Delta_2=v_2^2-4w_2=16(\lambda_1^2-\lambda_2^2+\lambda_3^2-\lambda_4^2)^2\left(S_2a^4+T_2a^2b^2+U_2b^4\right),$$
where
$$S_2=(\lambda_2^2-\lambda_4^2)^2>0, \qquad U_2=(\lambda_1^2-\lambda_3^2)^2>0,$$ 
$$T_2=2\left(-\,{\lambda_{{1}}}^{2}{\lambda_{{2}}}^{2}-\,{\lambda_{{1}}}^{2}{
\lambda_{{4}}}^{2}+2\,{\lambda_{{1}}}^{2}{\lambda_{{3}}}^{2}+2\,{
\lambda_{{2}}}^{2}{\lambda_{{4}}}^{2}-\,{\lambda_{{3}}}^{2}{\lambda_{
{4}}}^{2}-\,{\lambda_{{2}}}^{2}{\lambda_{{3}}}^{2}\right).$$
Moreover, the discriminant of the quadratic expression
$S_2a^4+T_2a^2b^2+U_2b^4$ is 
$$
T_2^2-4S_2U_2=-16(\lambda_1^2-\lambda_2^2)(\lambda_1^2-\lambda_4^2)(\lambda_2^2-\lambda_3^2)(\lambda_3^2-\lambda_4^2)<0
$$
which implies that $S_2a^4+T_2a^2b^2+U_2b^4>0$ and consequently $\Delta_2>0$. Therefore, the equation $t^2+v_2t+w_2=0$ has two non-zero distinct real roots of opposite signs and thus the equation $z^4+v_2z^2+w_2=0$ has two distinct real roots and two distinct purely imaginary roots. Thus $M_{a,b}^2 $ is a non-degenerate equilibrium and $\operatorname{span}_ \mathbb{R}\left\{\mathbf{D}X_{H|_{\operatorname{Orb}_{c_1;c_2}}}(M_{a,b}^2),\,  \mathbf{D}X_{I|_{\operatorname{Orb}_{c_1;c_2}}}(M_{a,b}^2) \right\}$ is a Cartan subalgebra of the second type in \eqref{cartan_subalgebras}. Thus, $M_{a,b}^2 $ is an unstable equilibrium of center-saddle type.

As before, the above computations being independent of the sign and permutation of $a$ and $b$, by an analogous reasoning we obtain non-degeneracy and instability for the other three equilibria $M_{-a,-b}^2,M_{b,a}^2,M_{-b,-a}^2$ in the Weyl orbit of $M_{a,b}^2$. 

\begin{thm}
All four equilibria in $\mathfrak{t}_2 \cap \operatorname{Orb}_{c_1; c_2}$ are non-degenerate, of center-saddle type and therefore unstable on the corresponding adjoint orbit. Thus these equilibria are also unstable for the Lie-Poisson dynamics on $\mathfrak{so}(4)$.
\end{thm}

\paragraph{The equilibria in $\mathfrak{t}_3$.} We proceed using the same techniques as in the previous cases. It is easy to see that $\mathbf{D}X_{I|_{\operatorname{Orb}_{c_1;c_2}}}(M_{a,b}^3)$ and $\mathbf{D}X_{H|_{\operatorname{Orb}_{c_1;c_2}}}(M_{a,b}^3)$ generate a 2-dimensional subspace. The eigenvalues of $\mathbf{D}X_{I|_{\operatorname{Orb}_{c_1;c_2}}}(M_{a,b}^3)$ are the roots of the equation
$$z^4+v_3z^2+w_3=0,$$
where
\begin{align*}
w_3&=16(\lambda_1^2-\lambda_3^2)(\lambda_1^2-\lambda_4^2)(\lambda_2^2-\lambda_3^2)(\lambda_2^2-\lambda_4^2)(a^2-b^2)^2>0 \\
v_3&=4\left[ (\lambda_1^2-\lambda_3^2)(\lambda_2^2-\lambda_3^2)+(\lambda_1^2-\lambda_4^2)(\lambda_2^2-\lambda_4^2) \right] 
{a}^{2}+\\
&\quad +4\left[ (\lambda_1^2-\lambda_3^2)(\lambda_1^2-\lambda_4^2)+(\lambda_2^2-\lambda_3^2)(\lambda_2^2-\lambda_4^2) \right] 
{b}^{2}>0.
\end{align*}
The quadratic equation $t^2+v_3t+w_3=0$ has the discriminant
$$\Delta_3=v_3^2-4w_3=16(\lambda_1^2-\lambda_3^2+\lambda_2^2-\lambda_4^2)^2\left(S_3a^4+T_3a^2b^2+U_3b^4\right),$$
where
$$S_3=(\lambda_3^2-\lambda_4^2)^2>0, \qquad U_3=(\lambda_1^2-\lambda_2^2)^2>0,$$ 
$$T_3=2\left[(\lambda_1^2-\lambda_3^2)(\lambda_2^2-\lambda_4^2)+(\lambda_1^2-\lambda_4^2)(\lambda_2^2-\lambda_3^2)\right]>0.$$
which implies that $S_3a^4+T_3a^2b^2+U_3b^4>0$ and consequently $\Delta_3>0$. 

Since $v_3>0,w_3>0$, the equation $t^2+v_3t+w_3=0$ has two non-zero distinct negative real roots and therefore equation $z^4+v_3z^2+w_3=0$ has two distinct pairs of purely imaginary roots and zero is not a root.
Thus  $\operatorname{span}_ \mathbb{R}\left\{\mathbf{D}X_{H|_{\operatorname{Orb}_{c_1;c_2}}}(M^3_{a,b}),\,
\mathbf{D}X_{I|_{\operatorname{Orb}_{c_1;c_2}}}(M^3_{a,b}) \right\}$ is a Cartan subalgebra an it is of the first type in \eqref{cartan_subalgebras}. It follows that the equilibrium $M_{a,b}^3$ is non-degenerate and nonlinearly stable because it is of center-center type (see \cite{BoFo04}, Theorem 1.5). 

The same holds for the other three equilibria $M_{-a,-b}^3$, $M_{b,a}^3$, $M_{-b,-a}^3$ in the Weyl group orbit of $M_{a,b}^3$. 

\begin{thm}
All four equilibria in $\mathfrak{t}_3 \cap \operatorname{Orb}_{c_1; c_2}$ are non-degenerate of center-center type and therefore nonlinearly stable on the corresponding adjoint orbit. These equilibria are also nonlinearly stable
for the Lie-Poisson dynamics on $\mathfrak{so}(4)$.
\end{thm}

Next, we begin the analysis of the remaining equilibria.

\paragraph{The equilibria in $\mathfrak{s} _\pm$.} 
The equilibria  from the families $\mathfrak{s}_+$ and $\mathfrak{s}_-$ are not isolated on the adjoint orbits. In fact, they come in curves or points described by intersecting the ellipsoids \eqref{first_condition_b} with the spheres \eqref{second_condition_b}, both families having the center at the origin. 

The characteristic 
equation of the linearized system at an equilibrium $M_e \in \operatorname{Orb}_{c_1;c_2}\cap 
\left[(\mathfrak{s}_+ \cup \mathfrak{s}_-) \setminus 
(\mathfrak{t}_1 \cup\mathfrak{t}_2 \cup \mathfrak{t}_3) 
\right]$ is
$$t^4(k_4t^2+k_1a_1^2+k_2a_2^2+k_3a_3^2)=0,$$
where
\begin{align*}
k_4&=\left( \lambda_{{2}}+\lambda_{{3}} \right) ^{2} \left( \lambda_{{1}}+\lambda_{{4}} \right) ^{2} \left( \lambda_{{3}}+\lambda_{{4}} \right) 
^{2} \left( \lambda_{{2}}+\lambda_{{4}} \right) ^{2} \left( \lambda_{{
3}}+\lambda_{{1}} \right) ^{2} \left( \lambda_{{1}}+\lambda_{{2}}
 \right) ^{2} >0 \\
k_1&=4\left( \lambda_{{3}}+\lambda_{{4}} \right)  \left( \lambda_{{2}}+
\lambda_{{4}} \right)  \left( \lambda_{{3}}+\lambda_{{1}} \right) 
 \left( \lambda_{{1}}+\lambda_{{2}} \right)  \left( \lambda_{{1}}
\lambda_{{4}}-\lambda_{{2}}\lambda_{{3}}\right)^{2}\geq 0 \\
k_2&=4\left( \lambda_{{3}}+\lambda_{{4}} \right)  
\left( \lambda_{{2}}+\lambda_{{3}} \right)  
\left( \lambda_{{1}}+\lambda_{{4}} \right) 
 \left( \lambda_{{1}}+\lambda_{{2}} \right)  
 \left( \lambda_{{3}}\lambda_{{1}}-\lambda_{{2}}\lambda_{{4}} \right) ^{2} > 0 \\
k_3&=4\left( \lambda_{{2}}+\lambda_{{4}}\right)  
\left( \lambda_{{2}}+\lambda_{{3}} \right)  
\left( \lambda_{{1}}+\lambda_{{4}} \right) 
\left(\lambda_{{3}}+\lambda_{{1}}\right)\left(\lambda_{{2}}
\lambda_{{1}}-\lambda_{{3}}\lambda_{{4}} \right)^{2} > 0
\end{align*}
and thus there are four zero eigenvalues; recall $\lambda_1>
\lambda_2 > \lambda_3 > \lambda_4 $. A double zero eigenvalue is expected since the generic orbit is four
dimensional. Restricting the linearized system to the 
tangent space to the orbit (which equals $\ker \mathbf{d}C_1
(M) \cap \ker \mathbf{d}C_2(M)$) yields a linear system
whose eigenvalues are the roots of the polynomial
$t^2(k_4t^2+k_1a_1^2+k_2a_2^2+k_3a_3^2)=0$. 
Therefore, the linearization
of the integrable system \eqref{sored} on the four 
dimensional adjoint orbit $\operatorname{Orb}_{c_1;c_2}$ 
at an equilibrium $M_e \in \operatorname{Orb}_{c_1;c_2}\cap 
\left[(\mathfrak{s}_+ \cup \mathfrak{s}_-) \setminus 
(\mathfrak{t}_1 \cup\mathfrak{t}_2 \cup \mathfrak{t}_3) 
\right]$ has the following eigenvalues: $0$ is a double 
eigenvalue and there are two other purely imaginary 
conjugate eigenvalues which can also degenerate to $0$.  Consequently, these equilibria can
only be degenerate cases of type 1 or type 2 in 
\eqref{cartan_subalgebras}. Thus, we cannot infer any
stability conclusion from the linearized system.

Note that the only time that $0$ can be a quadruple
eigenvalue is when $a_2=a_3=k_1=0$.

As before, we use the additional 
constant of motion $I_{\operatorname{Orb}_{c_1;c_2}}:=
I|_{\operatorname{Orb}_{c_1;c_2}}$ that commutes with
$H_{\operatorname{Orb}_{c_1;c_2}}: = 
H|_{\operatorname{Orb}_{c_1;c_2}}$. However, by 
Proposition \ref{two_types_equilibria}, 
$\mathbf{d}I_{\operatorname{Orb}_{c_1;c_2}}(M_e)\neq 0$, so we can not apply the method used for studying the stability for the equilibria in $\mathcal{K}_0 = \operatorname{Orb}_{c_1;c_2} \cap (\mathfrak{t}_1 \cup \mathfrak{t}_2 \cup\mathfrak{t}_3)$. 

We shall use energy methods (see \cite{Arnold65}, 
\cite{HoMaRaWe85}, \cite{OrRa}, \cite{BiPu07}). 
If
$$m_0 = -\frac{1}{\lambda_1+\lambda_2+\lambda_3+\lambda_4}\,, \qquad 
n_0 = -\frac{1}{\lambda_1+\lambda_2+\lambda_3+\lambda_4}$$
then $\mathbf{d}(H+m_0C_1+n_0C_2)(M_e)=0$ and 
the Hessian 
$\mathbf{D}^2(H+m_0C_1+n_0C_2)(M_e) $ has 
characteristic
polynomial 
$$t^3(t-\alpha_1)(t-\alpha_2)(t-\alpha_3)=0,$$
where
\begin{align*}
\alpha_1&=\frac{(\lambda_1+\lambda_3)^2+(\lambda_2+\lambda_4)^2}{(\lambda_1+\lambda_3)(\lambda_2+\lambda_4)(\lambda_1+\lambda_2+\lambda_3+\lambda_4)} >0\\
\alpha_2&=\frac{(\lambda_1+\lambda_4)^2+(\lambda_2+\lambda_3)^2}{(\lambda_1+\lambda_4)(\lambda_2+\lambda_3)(\lambda_1+\lambda_2+\lambda_3+\lambda_4)}>0\\
\alpha_3&=\frac{(\lambda_1+\lambda_2)^2+(\lambda_3+\lambda_4)^2}{(\lambda_1+\lambda_2)(\lambda_3+\lambda_4)(\lambda_1+\lambda_2+\lambda_3+\lambda_4)}>0.
\end{align*}
We suppose, without lose of generality, that $a_1\not= 0,a_2\not= 0$.
The Hessian $\mathbf{D}^2H_{\operatorname{Orb}_
{c_1;c_2}}(M_e)=\mathbf{D}^2(H+m_0C_1+n_0C_2)(M_e)|_{T_{M_e}
\operatorname{Orb}_{c_1;c_2}}$,  where $T_{M_e}\operatorname
{Orb}_{c_1;c_2}=\ker \mathbf{d}C_1(M_e)\cap \ker \mathbf{d}
C_2(M_e)$ has eigenvalues $0,\beta_1,\beta_2,\beta_3$ computed in a conveniently chosen basis for $T_{M_e}\operatorname
{Orb}_{c_1;c_2}$. Using Vi\`ete's relations for the characteristic polynomial of $\mathbf{D}^2H_{\operatorname{Orb}_{c_1;c_2}}(M_e)$ 
computed in the above basis we have:
$$\beta_1\beta_2\beta_3=4\frac{A_1+A_2\left(\frac{a_3}{a_1}\right)^2+A_3\left(\frac{a_3}{a_2}\right)^2}{B},$$
where:
\begin{align*}
A_1&=(\lambda_1^2-\lambda_2^2)^2(\lambda_3^2-\lambda_4^2)^2\left[(\lambda_1+\lambda_2)^2+(\lambda_3+\lambda_4)^2\right]>0;\\
A_2&=(\lambda_1^2-\lambda_4^2)^2(\lambda_2^2-\lambda_3^2)^2(\lambda_1+\lambda_4)^2>0;\\
A_3&=(\lambda_1^2-\lambda_3^2)^2(\lambda_2^2-\lambda_4^2)^2(\lambda_2+\lambda_4)^2>0;
\end{align*}
$$B=(\lambda_1+\lambda_2)^3(\lambda_1+\lambda_3)(\lambda_1+\lambda_4)(\lambda_2+\lambda_3)(\lambda_2+\lambda_4)(\lambda_3+\lambda_4)^3\cdot$$
$$\cdot (\lambda_1-\lambda_2)^2(\lambda_1+\lambda_2+\lambda_3+\lambda_4)>0.$$
This shows that $\beta_1,\beta_2,\beta_3$ are all non-zero and positive since $\mathbf{D}^2H_{\operatorname{Orb}_{c_1;c_2}}(M_e)$ is positive semi-definite as it is a restriction of the positive semi-definite bilinear form $\mathbf{D}^2(H+m_0C_1+n_0C_2)(M_e) $.

Now we can conclude the following theorem.

\begin{thm}\label{second_stability_result}
Each curve of equilibria in $\operatorname{Orb}_{c_1;c_2}\cap 
\left[(\mathfrak{s}_+ \cup \mathfrak{s}_-) \setminus 
(\mathfrak{t}_1 \cup\mathfrak{t}_2 \cup \mathfrak{t}_3) 
\right]$ is nonlinearly stable. That is, if a solution
of \eqref{sored} starts near an equilibrium on such a curve,
at any later time it will stay close to the curve in 
$\operatorname{Orb}_{c_1;c_2}\cap 
\left[(\mathfrak{s}_+ \cup \mathfrak{s}_-) \setminus 
(\mathfrak{t}_1 \cup\mathfrak{t}_2 \cup \mathfrak{t}_3) 
\right]$ containing this equilibrium, but in the direction
of this curve it may drift. Since 
$\operatorname{Orb}_{c_1;c_2}$ is a generic adjoint 
orbit, if the perturbation is close to the given 
equilibrium but on a neighboring adjoint orbit the same 
situation occurs.
\end{thm}

\begin{rem}
One can pose the legitimate question if the statement of
the theorem above could be strengthened in the sense that
the drift in the neutral direction is impossible, at least
for some equilibria. This would
then prove the nonlinear stability of such an 
equilibrium on these curves of equilibria. To achieve this,
one would have to show that $\mathbf{D}^2H|_{\operatorname{Orb}_{c_1;c_2}}(M_e)$ is definite when restricted to the leaf $L_{M_e}$ (to be defined below), which would give nonlinear stability by Arnold's method (which is proved to be equivalent with the other energy methods, see \cite{BiPu07}). We shall show below that the method is
inconclusive so we do not know which, if any, of the equilibria on these curves are nonlinearly stable.

So let's try to apply the Arnold stability method to such
an equilibrium $M_e$. We need to study the definiteness of the Hessian of the constant of the motion $H_{\operatorname{Orb}_{c_1;c_2}}+\alpha I_{\operatorname{Orb}_{c_1;c_2}}$
evaluated at $M_e $ restricted to the invariant level set $L_{M_e}:=I_{\operatorname{Orb}_{c_1;c_2}}^{-1}\left(I_{\operatorname{Orb}_{c_1;c_2}}(M_e)\right)$ of the dynamics \eqref{sored}.  The conditions $\mathbf{d} (H_{\operatorname{Orb}_{c_1;c_2}}+\alpha I_{\operatorname{Orb}_{c_1;c_2}})(M_e)=0$, $\mathbf{d} H_{\operatorname{Orb}_{c_1;c_2}}(M_e)= 0$, and $\mathbf{d} I_{\operatorname{Orb}_{c_1;c_2}}(M_e)\neq 0$ (since $M_e \in \mathcal{K}_1$) imply $\alpha =0$. Consequently, we have to study the definiteness of 
the Hessian of $H_{\operatorname{Orb}_{c_1;c_2}}$ at $M_e $ restricted to the tangent space at $M_e$ of $L_{M_e}$. 
We shall prove below that this definiteness does not hold. 

Let $c_{M_e}(t)$ be the curve of equilibria for $X_{H_{\operatorname{Orb}_{c_1;c_2}}}$ with $c_{M_e}(0)=M_e$. Then $X_{H_{\operatorname{Orb}_{c_1;c_2}}}(c_{M_e}(t))=0$ and by differentiation  $\mathbf{D}X_{H_{\operatorname{Orb}_{c_1;c_2}}}(0)\cdot \dot c_{M_e}(0)=0$. Equivalently, using the formula of the linearization of a Hamiltonian vector field at a critical point on a symplectic manifold, we have $\Lambda^{-1} \mathbf{D}^2H|_{\operatorname{Orb}_{c_1;c_2}}(M_e)\cdot \dot c_{M_e}(0)=0$, where $\Lambda$ is the $4\times 4$ matrix associated to the symplectic form on the adjoint orbit $\operatorname{Orb}_{c_1;c_2}$. As $\Lambda$ is nondegenerate, we obtain that $\mathbf{D}^2H|_{\operatorname{Orb}_{c_1;c_2}}(M_e)\cdot \dot c_{M_e}(0)=0$, which shows that $\dot c_{M_e}(0)$ is in the eigendirection corresponding to the eigenvalue $0$ for $\mathbf{D}^2H|_{\operatorname{Orb}_{c_1;c_2}}(M_e)$.

We shall prove that $X_{I_{\operatorname{Orb}_{c_1;c_2}}}(M_e)$ is collinear with $\dot c_{M_e}(0)$. Suppose not; then $X_{I_{\operatorname{Orb}_{c_1;c_2}}}(M_e)$ is not tangent to the curve $c_{M_e}(t)$ at the point $M_e$. Consequently, for a small $s\in \mathbf{R}$, along the integral curve of $X_{I_{\operatorname{Orb}_{c_1;c_2}}}$, 
we find an $s \in \mathbb{R}$ such that  
$\Phi_s^{I_{\operatorname{Orb}_{c_1;c_2}}}(M_e)=x_1$, where $x_1\not\in c_{M_e}(t)$. Since $\left\{H_{\operatorname{Orb}_{c_1;c_2}}, I_{\operatorname{Orb}_{c_1;c_2}}\right\} = 0$, the flows $\Phi_t^{H_{\operatorname{Orb}_{c_1;c_2}}} $ of $X_{H_{\operatorname{Orb}_{c_1;c_2}}}$ and $\Phi_s^{I_{\operatorname{Orb}_{c_1;c_2}}}$ of $X_{I_{\operatorname{Orb}_{c_1;c_2}}}$ commute and hence 
$$\Phi_{-s}^{I_{\operatorname{Orb}_{c_1;c_2}}}\circ \Phi_{-t}^{H_{\operatorname{Orb}_{c_1;c_2}}}\circ \Phi_s^{I_{\operatorname{Orb}_{c_1;c_2}}}\circ \Phi_t^{H_{\operatorname{Orb}_{c_1;c_2}}}(M_e)=M_e.$$
Because in a neighborhood of $M_e$ the only equilibria for the vector field $X_{H_{\operatorname{Orb}_{c_1;c_2}}}$ are of on the curve $c_{M_e}(t)$, we conclude that $x_2:=\Phi_{-t}^{H_{\operatorname{Orb}_{c_1;c_2}}}(x_1)\neq x_1$. But then
\begin{align*}
\Phi_{-s}^{I_{\operatorname{Orb}_{c_1;c_2}}}(x_2)
&=\Phi_{-s}^{I_{\operatorname{Orb}_{c_1;c_2}}}\left(\Phi_{-t}^{H_{\operatorname{Orb}_{c_1;c_2}}}(x_1)\right)
= \left(\Phi_{-s}^{I_{\operatorname{Orb}_{c_1;c_2}}}\circ 
\Phi_{-t}^{H_{\operatorname{Orb}_{c_1;c_2}}}\right)
\left(\Phi_s^{I_{\operatorname{Orb}_{c_1;c_2}}}(M_e) \right)\\
&=\left(\Phi_{-s}^{I_{\operatorname{Orb}_{c_1;c_2}}}\circ \Phi_{-t}^{H_{\operatorname{Orb}_{c_1;c_2}}}\circ \Phi_s^{I_{\operatorname{Orb}_{c_1;c_2}}}\circ \Phi_t^{H_{\operatorname{Orb}_{c_1;c_2}}}\right)(M_e)=M_e
\end{align*} 
since $\Phi_t^{H_{\operatorname{Orb}_{c_1;c_2}}}(M_e)=M_e$, 
because $M_e$ is an equlibrium of \eqref{sored}.
Therefore, $\Phi_s^{I_{\operatorname{Orb}_{c_1;c_2}}}(M_e)=x_2$. However, we also have $\Phi_s^{I_{\operatorname{Orb}_{c_1;c_2}}}(M_e)=x_1$, which contradicts the uniqueness of
integral curves of $X_{I_{\operatorname{Orb}_{c_1;c_2}}}$. Therefore $X_{I_{\operatorname{Orb}_{c_1;c_2}}}(M_e)$ is collinear with $\dot{c}_{M_e}(0)$ and is thus in 
the eigendirection corresponding to the eigenvalue $0$ for 
$\mathbf{D}^2H|_{\operatorname{Orb}_{c_1;c_2}}(M_e)$.

We have $T_{M_e}L_{M_e}\oplus \mathbf{D} I_{\operatorname{Orb}_{c_1;c_2}}(M_e)=T_{M_e}\operatorname{Orb}_{c_1;c_2}$ and, using Darboux coordinates, it is easily shown that we also have $X_{I_{\operatorname{Orb}_{c_1;c_2}}}\perp \mathbf{D} I_{\operatorname{Orb}_{c_1;c_2}}$, which shows that 
$X_{I_{\operatorname{Orb}_{c_1;c_2}}}(M_e)\in T_{M_e}L_{M_e}$. 

Thus the Hessian of $H_{\operatorname{Orb}_{c_1;c_2}}$ restricted to the level manifold $L_{M_e}$ has a $0$ eigenvalue in the direction $X_{I_{\operatorname{Orb}_{c_1;c_2}}}(M_e)$. 

This shows that the use of the constant of the motion $I$ does not improve the nonlinear stability result of equilibrium points of type $\mathcal{K}_1$ in 
Theorem \ref{second_stability_result}.
\end{rem}

\begin{rem}
The eigenvalue $0$ for the Hessian is expected since
the equilibrium $M_e$ lies on the curve obtained by 
intersecting the ellipsoids \eqref{second_condition} and 
\eqref{first_condition}, both having centers at the origin; 
the $0$-eigenspace is tangent to this curve of equilibria, as proved above.
So the only stability we can expect is the stability
transversal to the direction of the curve of equilibria.  
\end{rem}

\noindent\textbf{Acknowledgements.} 
The authors acknowledge the very helpful discussions with Daisuke Tarama.
T.S. Ratiu was partially supported by Swiss NSF grant 200020-126630 and by the government grant of the Russian Federation for support of  research projects implemented by leading scientists, Lomonosov Moscow  State University, under the agreement No. 11.G34.31.0054.
P. Birtea has been supported by CNCSIS - UEFISCDI, project number PNII
- IDEI code 1081/2008 No. 550/2009 and by Romanian National University Research Council (CNCSIS-UEFISCU), project PNII-IDEI 131/2008.
M. Turhan has been partially supported by a Swiss National Science Foundation grant.\\
Some of the computations 
have been done using Maple 11.

\end{document}